# Characterizations of Mixed Herz-Hardy Spaces and their Applications

Yichun Zhao, Mingquan Wei, Jiang Zhou *

**Abstract:** The purpose of this paper is to introduce and investigate some basic properties of mixed homogeneous Herz-Hardy spaces $H\dot{K}_{\vec{p}}^{\alpha,q}(\mathbb{R}^n)$ and mixed non-homogeneous Herz-Hardy spaces $HK_{\vec{p}}^{\alpha,q}(\mathbb{R}^n)$. Furthermore, we establish the atom and molecular decompositions for $H\dot{K}_{\vec{p}}^{\alpha,q}(\mathbb{R}^n)$ and $HK_{\vec{p}}^{\alpha,q}(\mathbb{R}^n)$, by which the boundedness for a wide class of sublinear operators on mixed Herz-Hardy spaces is obtained. As a byproduct, the dual spaces of mixed homogeneous Herz-Hardy spaces are deduced.

## 1 Introduction

The study of Herz spaces originated from the work of Beurling [1]. Later, Herz spaces were systematically studied by Herz [2] to study the Fourier series and Fourier transform. In the 1990s, Lu and Yang [3,4] introduced the homogeneous Herz spaces $\dot{K}_q^{\alpha,p}(\mathbb{R}^n)$ and non-homogeneous Herz spaces $K_q^{\alpha,p}(\mathbb{R}^n)$ and extended the boundedness of a large class of operators to these spaces.

As is well known, Hardy spaces are proper substitutes of Lebesgue spaces in some situations. For example, the Riesz transforms are bounded on the Hardy space $H^p(\mathbb{R}^n)$, but not bounded on the corresponding Lebesgue space $L^p(\mathbb{R}^n)$ when $0 < p \leqslant 1$. The theory of classical Hardy spaces was originally studied by Stein and Weiss [5] and then systematically developed in [6]. Hardy spaces also have various real variables characterizations, such as different maximal operators characterizations, atom characterization [7] and molecular characterization [8], and so on. These characterizations greatly facilitate the researchers to derive the dual spaces of Hardy spaces and establish the boundedness of operators on these spaces, see [9] for more details. As a variant of the classical real Hardy spaces, the Hardy spaces associated with the Beurling algebras on the real line were first introduced by Chen and Lau [10], in which the dual spaces and the maximal function characterizations of these spaces were established. Later, García-Cuerva [11] extended the results of Chen and Lau [10] to higher-dimensional case, and García-Cuerva and Herrero [12] further studied the maximal function, atom, and Littlewood-Paley function characterizations of these Herz-Hardy type spaces. In 1995, Lu and Yang [13] systematically studied Herz-Hardy spaces with general indices and established their atom and molecular characterizations. For more studies of Herz-Hardy type spaces, the readers can refer to [14–18].

---
*Corresponding author. The research was supported by Natural Science Foundation of China (Grant No. 12061069) and the Natural Science Foundation of Henan Province(Grant No. 202300410338 ).

Recently, mixed-norm Lebesgue spaces, as natural extensions of the classical Lebesgue spaces, have attracted widespread attention. The theory of mixed-norm function spaces can be traced back to the work of Benedek and Panzone [19], in which it was proved that $L^{\vec{p}}(\mathbb{R}^n)$ also possesses some basic properties similar to the classical Lebesgue spaces, such as completeness, Hölder's inequality, Minkowski's inequality, and so on. These properties provide the possibility to solve a series of subsequence problems. On the other hand, since the study of partial differential equations (for example the heat equation and the wave equation) always involves both space and time variables, mixed-norm spaces possess better structures than classical spaces in the time-space estimates for PDEs. For these reasons, many researchers renewed the interest in mixed-norm Lebesgue spaces and extended them to other mixed-norm function spaces. For instance, the real variable characterizations and the atom characterization of mixed-norm Hardy spaces were established in [20, 21] In 2019, Nogayama [22] introduced mixed-norm Morrey spaces and gave some applications in the operator theory.

Note that recently Herz spaces were also extended to the mixed-norm situation by Wei [23]. By extending the extrapolation theory to mixed Herz spaces, Wei [23] established the boundedness of some classical operators in harmonic analysis on these spaces. Moreover, the extrapolation theory can further give boundedness results of some classical operators on mixed Herz spaces.

Motivated by [23] and the theory of Herz-type Hardy spaces, we introduce mixed Herz-Hardy spaces and investigate some basic properties of these spaces in this paper. Moreover, we will establish the atom decomposition and molecular decomposition for these spaces. As applications, the boundedness for a wide class of sublinear operators on mixed Herz-Hardy spaces is obtained, and in addition, the dual spaces of mixed homogeneous Herz-Hardy spaces are also deduced.

The organization of the remainder of this article is as follows. Some necessary definitions and lemmas are given in Section 2. The atom decomposition and the molecular decomposition of mixed Herz-Hardy spaces will be given in Section 3 and Section 4, respectively. Some applications of the atom and molecular decomposition are presented in Section 5.

## 2 Preliminary

Throughout this paper, we use the following notations. The letter $\vec{q}$ will denote $n$-tuples of the numbers in $(0, \infty]$ $(n \geqslant 1)$, $\vec{q} = (q_1, q_2, \ldots, q_n)$. By definition, the inequality $0 < \vec{q} < \infty$ means that $0 < q_i < \infty$ for all $i$. For $\vec{q} = (q_1, q_2, \ldots, q_n)$, write

$$\frac{1}{\vec{q}} = \left(\frac{1}{q_1}, \frac{1}{q_2}, \ldots, \frac{1}{q_n}\right), \quad \vec{q'} = (q'_1, q'_2, \ldots, q'_n),$$

where $q'_i = q_i/(q_i - 1)$ is conjugate exponent of $q_i$. $|B|$ denotes the volume of the ball $B$, $\chi_E$ is the characteristic function of a set $E$. $A \sim B$ means that $A \lesssim B$ and $B \lesssim A$, $[a]$ denotes take the integer number for $a$. Let $B_k = \{x \in \mathbb{R}^n : |x| \leqslant 2^k\}$ and $A_k = B_k \backslash B_{k-1}$ for any $k \in \mathbb{Z}$. Denote $\chi_k = \chi_{A_k}$ for any $k \in \mathbb{Z}$, and $\widetilde{\chi}_k = \chi_k$ for any $k \in \mathbb{N}$, $\widetilde{\chi}_0 = \chi_{B_0}$.

**Definition 2.1.** *(Mixed Lebesgue spaces)( [19]) Let $\vec{p} = (p_1, p_2, \ldots, p_n) \in (0, \infty]^n$. Then the mixed Lebesgue space $L^{\vec{p}}(\mathbb{R}^n)$ is defined to be the set of all measurable functions $f$ such that*

$$\|f\|_{L^{\vec{p}}(\mathbb{R}^n)} := \left(\int_{\mathbb{R}} \cdots \left(\int_{\mathbb{R}} \left(\int_{\mathbb{R}} |f(x_1, x_2, \ldots, x_n)|^{p_1} dx_1\right)^{\frac{p_2}{p_1}} dx_2\right)^{\frac{p_3}{p_2}} \ldots dx_n\right)^{\frac{1}{p_n}} < \infty,$$

*If $p_j = \infty$, then we have to make appropriate modifications.*

**Definition 2.2.** ( [23]) *Let $\alpha \in \mathbb{R}$, $0 < q \leqslant \infty$, $0 < \vec{p} \leqslant \infty$. The mixed homogeneous Herz space $\dot{K}_{\vec{p}}^{\alpha,q}(\mathbb{R}^n)$ is defined by*

$$\dot{K}_{\vec{p}}^{\alpha,q}(\mathbb{R}^n) := \left\{ f \in L_{\text{loc}}^{\vec{p}}(\mathbb{R}^n) : \|f\|_{\dot{K}_{\vec{p}}^{\alpha,q}(\mathbb{R}^n)} = \left( \sum_{k \in \mathbb{Z}} 2^{k\alpha q} \|f\chi_k\|_{L^{\vec{p}}(\mathbb{R}^n)}^q \right)^{1/q} < \infty \right\}.$$

**Definition 2.3.** ( [23]) *Let $\alpha \in \mathbb{R}$, $0 < q \leqslant \infty$, $0 < \vec{p} \leqslant \infty$. The non-homogeneous mixed Herz space $K_{\vec{p}}^{\alpha,q}(\mathbb{R}^n)$ is defined by*

$$K_{\vec{p}}^{\alpha,q}(\mathbb{R}^n) := \left\{ f \in L_{\text{loc}}^{\vec{p}}(\mathbb{R}^n) : \|f\|_{K_{\vec{p}}^{\alpha,q}(\mathbb{R}^n)} = \left( \sum_{k=0}^{\infty} 2^{k\alpha q} \|f\widetilde{\chi}_k\|_{L^{\vec{p}}(\mathbb{R}^n)}^q \right)^{1/q} < \infty \right\}.$$

**Remark 2.1.** (i) *If $0 < \vec{p} = (p_1, p_2, \ldots, p_n) \leqslant \infty$ and $p_1 = p_2 = \cdots = p_n = p$, then $\dot{K}_{\vec{p}}^{\alpha,q}(\mathbb{R}^n) = \dot{K}_p^{\alpha,q}(\mathbb{R}^n)$ or $K_{\vec{p}}^{\alpha,q}(\mathbb{R}^n) = K_p^{\alpha,q}(\mathbb{R}^n)$, where $\dot{K}_p^{\alpha,q}(\mathbb{R}^n)$ and $K_p^{\alpha,q}(\mathbb{R}^n)$ are classical Herz spaces.*

(ii) *The mixed homogeneous Herz space $\dot{K}_{\vec{p}}^{\alpha,q}(\mathbb{R}^n)$ and the non-homogeneous mixed Herz space $K_{\vec{p}}^{\alpha,q}(\mathbb{R}^n)$ are quasi-Banach spaces. But, if $\vec{p}, q \geqslant 1$, they are Banach spaces. These can be inferred from definitions of mixed Lebesgue spaces and classical Herz spaces.*

Now we give the definition of mixed Herz-Hardy spaces $H\dot{K}_{\vec{p}}^{\alpha,q}(\mathbb{R}^n)$ and $HK_{\vec{p}}^{\alpha,q}(\mathbb{R}^n)$. Denote $\mathcal{S}(\mathbb{R}^n)$ by the Schwartz space of all rapidly decreasing infinitely differentiable functions on $\mathbb{R}^n$, and $\mathcal{S}'(\mathbb{R}^n)$ by the dual space of $\mathcal{S}(\mathbb{R}^n)$.

Let $\mathscr{M}_N f(x)$ be the grand maximal function of $f(x)$ defined by

$$\mathscr{M}_N f(x) = \sup_{\phi \in \mathcal{A}_N} M_1^*(f;\phi)(x),$$

where $\mathcal{A}_N = \left\{ \phi \in \mathcal{S}(\mathbb{R}^n) : \sup_{|\alpha|,|\beta| \leqslant N} |x^\alpha D^\beta \phi(x)| \leqslant 1 \right\}$ and $N > n+1$, $M_1^*(f;\phi)$ is a special nontangential maximal operator defined by

$$M_1^*(f;\phi)(x) = \sup_{|y-x|<t} |\phi_t * f(y)|$$

with $\phi_t(x) = t^{-n}\phi(x/t)$.

**Definition 2.4.** *Let $\alpha \in \mathbb{R}, 0 < q < \infty, 1 < \vec{p} \leqslant \infty$, and $N > n+1$. The mixed homogeneous Herz-type Hardy space $H\dot{K}_{\vec{p}}^{\alpha,p}(\mathbb{R}^n)$ is defined by*

$$H\dot{K}_{\vec{p}}^{\alpha,q}(\mathbb{R}^n) = \left\{ f \in \mathcal{S}'(\mathbb{R}^n) : \mathscr{M}_N f(x) \in \dot{K}_{\vec{p}}^{\alpha,q}(\mathbb{R}^n) \right\},$$

*that is*

$$\|f\|_{H\dot{K}_{\vec{p}}^{\alpha,q}(\mathbb{R}^n)} = \|\mathscr{M}_N f\|_{\dot{K}_{\vec{p}}^{\alpha,q}(\mathbb{R}^n)}.$$

**Definition 2.5.** *Let $\alpha \in \mathbb{R}, 0 < q < \infty, 1 < \vec{p} \leqslant \infty$, and $N > n+1$. The mixed non-homogeneous Herz-type Hardy space $HK_{\vec{p}}^{\alpha,p}(\mathbb{R}^n)$ is defined by*

$$HK_{\vec{p}}^{\alpha,q}(\mathbb{R}^n) = \left\{ f \in \mathcal{S}'(\mathbb{R}^n) : \mathscr{M}_N f(x) \in K_{\vec{p}}^{\alpha,q}(\mathbb{R}^n) \right\},$$

*that is*

$$\|f\|_{H\dot{K}_{\vec{p}}^{\alpha,q}(\mathbb{R}^n)} = \|\mathscr{M}_N f\|_{\dot{K}_{\vec{p}}^{\alpha,q}(\mathbb{R}^n)}.$$

**Remark 2.2.** (i) *Notice that if $q = 1, 1 < \vec{p} < \infty$, and $\alpha = n - \sum_{i=1}^{n} 1/p_i$, then the space $H\dot{K}_{\vec{p}}^{\alpha,1}(\mathbb{R}^n)$ is the space $HA^{\vec{p}}(\mathbb{R}^n)$ defined in [24], which can be inferred from the atom decomposition of mixed Herz-Hardy spaces (see section 3).*

(ii) *By the sublinearity of $\mathscr{M}_N$ and the definition of mixed Herz-Hardy spaces, we can conclude that mixed Herz-Hardy spaces are quasi-Banach spaces.*

Now, we will give some necessary lemmas and basic properties of mixed Herz-Hardy spaces.

Besides the grand maximal operator $\mathscr{M}_N$, the space $H\dot{K}_{\vec{p}}^{\alpha,q}(\mathbb{R}^n)$ can also be characterized by some other maximal-type operators. To give these characterizations, we first introduce some maximal-type operators. Let $\phi \in \mathcal{S}(\mathbb{R}^n)$ such that $\int_{\mathbb{R}^n} \phi(x)dx = 1$. For $t > 0$, $f \in \mathcal{S}'(\mathbb{R}^n)$, define the smooth maximal operator $M(f;\phi)$ by

$$M(f;\phi)(x) = \sup_{t>0} |(f \star \phi_t)(x)|.$$

Also, we define the non-tangential maximal operator $M_a^*(f;\phi)$ (with $a > 1$) and the auxiliary maximal operator $M_b^{**}(f;\phi)$ (with $M \in \mathbb{N}$) by

$$M_a^*(f;\phi)(x) = \sup_{t>0} \sup_{|x-y|<at} |f \ast \phi_t(x)|$$

and

$$M_b^{**}(f;\phi)(x) = \sup_{(y,t) \in \mathbb{R}_+^{n+1}} |f \ast \phi_t(y)| \left(\frac{t}{|x-y|+t}\right)^b.$$

**Proposition 2.1.** *Let $0 < \alpha < \infty, 0 < q < \infty$ and $1 < \vec{p} < \infty$. For $f \in \mathcal{S}'(\mathbb{R}^n)$, the following statements are equivalent:*

(i) $f \in H\dot{K}_{\vec{p}}^{\alpha,q}(\mathbb{R}^n)$.

(ii) *There exists a function $\phi \in \mathcal{S}(\mathbb{R}^n)$ and $\int_{\mathbb{R}^n} \phi(x)dx = 1$, such that for some $a > 1$, $M_a^*(f;\phi) \in \dot{K}_{\vec{p}}^{\alpha,q}(\mathbb{R}^n)$.*

(iii) *There exists a function $\phi \in \mathcal{S}(\mathbb{R}^n)$ and $\int_{\mathbb{R}^n} \phi(x)dx = 1$, such that $M_1^*(f;\phi) \in \dot{K}_{\vec{p}}^{\alpha,q}(\mathbb{R}^n)$.*

(iv) *There exists a function $\phi \in \mathcal{S}(\mathbb{R}^n)$ and $\int_{\mathbb{R}^n} \phi(x)dx = 1$, such that $M(f;\phi) \in \dot{K}_{\vec{p}}^{\alpha,q}(\mathbb{R}^n)$.*

*Proof.* This proof is similar to that of the characterizations of $H^p(\mathbb{R}^n)$ by different maximal-type functions. One can see ([9], Section 2.1.2) for the details. □

**Lemma 2.1.** *Let $0 < q \leqslant \infty$, $1 < \vec{p} < \infty$, and $-\sum_{i=1}^{n} 1/p_i < \alpha < n\left(1 - 1/n \sum_{i=1}^{n} 1/p_i\right)$. Suppose a sublinear operators $T$ satisfied that*

(i) *$T$ is bounded on $L^{\vec{p}}(\mathbb{R}^n)$;*

(ii) *for any $f \in L^1(\mathbb{R}^n)$ with compact support have*

$$|Tf(x)| \leqslant C \int_{\mathbb{R}} \frac{|f(x)|}{|x-y|^n}, \qquad x \notin \text{supp} f. \tag{1}$$

*Then $T$ is bounded on $\dot{K}_{\vec{p}}^{\alpha,q}(\mathbb{R}^n)$.*

*Proof.* This proof is similar to that of ([23], Theorem 4.2), and the only difference is that we should take the size condition into consideration. For simplification, we omit the details. □

**Lemma 2.2.** ( [20]) Let $\vec{p} \in (0,\infty)^n$ and $0 < \theta < \min(1, p_1, \ldots, p_n)$. Then for every $\alpha > 0, \beta > n/\theta, N > 2n + \beta + 2$ and $\phi \in \mathcal{S}(\mathbb{R}^n)$ with $\int_{\mathbb{R}^n} \phi(x)dx \neq 0$, for all $f \in \mathcal{S}'(\mathbb{R}^n)$, there holds

$$\|M(f;P)\|_{L^{\vec{p}}(\mathbb{R}^n)} \sim \|M(f;\phi)\|_{L^{\vec{p}}(\mathbb{R}^n)} \sim \|M_\alpha^*(f;\phi)\|_{L^{\vec{p}}(\mathbb{R}^n)} \sim \|M_\beta^{**}(f;\phi)\|_{L^{\vec{p}}(\mathbb{R}^n)} \sim \|\mathscr{M}_N f\|_{L^{\vec{p}}(\mathbb{R}^n)}.$$

**Proposition 2.2.** Let $1 < \vec{p} < \infty$, $0 < \theta < \min(1, p_1, \ldots, p_n)$, $\beta > n/\theta$, and $N > 2n + \beta + 2$. Then $\mathscr{M}_N$ is bounded on $L^{\vec{p}}(\mathbb{R}^n)$.

*Proof.* By Lemma 2.2, $\|M(f;P)\|_{L^{\vec{p}}(\mathbb{R}^n)} \sim \|\mathscr{M}_N f\|_{L^{\vec{p}}(\mathbb{R}^n)}$. Then, we only need to consider the boundedness of operators $M(f;P)$ on $L^{\vec{p}}(\mathbb{R}^n)$.

If $\psi \in L^1(\mathbb{R}^n)$ is the decreasing radial dominant function of the Poisson function, then

$$M(f;P)(x) \leqslant C\|\psi\|_{L^1} M f(x).$$

Since the Hardy-Littlewood maximal operator is bounded on $L^{\vec{p}}(\mathbb{R}^n)$ (see [25]), we have

$$\|\mathscr{M}_N f\|_{L^{\vec{p}}(\mathbb{R}^n)} \sim \|M(f;P)\|_{L^{\vec{p}}(\mathbb{R}^n)} \leqslant C\|\psi\|_{L^1}\|Mf\|_{L^{\vec{p}}(\mathbb{R}^n)} \leqslant C\|f\|_{L^{\vec{p}}(\mathbb{R}^n)}.$$

□

**Proposition 2.3.** Let $0 < q < \infty, 1 < \vec{p} = (p_1, p_2, \ldots, p_n) < \infty$ and $N > n+1$. If $-\sum_{i=1}^m 1/p_i < \alpha < n(1 - 1/n\sum_{i=1}^m 1/p_i)$. Then

$$H\dot{K}_{\vec{p}}^{\alpha,q}(\mathbb{R}^n) \bigcap L_{\text{loc}}^{\vec{p}}(\mathbb{R}^n\backslash\{0\}) = \dot{K}_{\vec{p}}^{\alpha,q}(\mathbb{R}^n)$$

and

$$HK_{\vec{p}}^{\alpha,q}(\mathbb{R}^n) \bigcap L_{\text{loc}}^{\vec{p}}(\mathbb{R}^n) = K_{\vec{p}}^{\alpha,q}(\mathbb{R}^n).$$

*Proof.* It suffices to prove homogeneous spaces.

By using the trivial inequality $|f(x)| \leqslant C\mathscr{M}_N(f)(x)$, we get

$$\|f\|_{\dot{K}_{\vec{p}}^{\alpha,q}(\mathbb{R}^n)} \lesssim \left(\sum_{k \in \mathbb{Z}} 2^{k\alpha q}\|\mathscr{M}_N(f)(x)\chi_k\|_{L^{\vec{p}}(\mathbb{R}^n)}^q\right)^q = \|\mathscr{M}_N(f)\|_{\dot{K}_{\vec{p}}^{\alpha,q}(\mathbb{R}^n)} = \|f\|_{H\dot{K}_{\vec{p}}^{\alpha,q}(\mathbb{R}^n)},$$

which yields $H\dot{K}_{\vec{p}}^{\alpha,q}(\mathbb{R}^n) \bigcap L_{\text{loc}}^{\vec{p}}(\mathbb{R}^n\backslash\{0\}) \subset \dot{K}_{\vec{p}}^{\alpha,q}(\mathbb{R}^n)$.

From the boundedness of sublinear operators on mixed Herz spaces $\dot{K}_{\vec{p}}^{\alpha,q}(\mathbb{R}^n)$, we can check that $\dot{K}_{\vec{p}}^{\alpha,q}(\mathbb{R}^n) \subset H\dot{K}_{\vec{p}}^{\alpha,q}(\mathbb{R}^n) \bigcap L_{\text{loc}}^{\vec{p}}(\mathbb{R}^n\backslash\{0\})$.

We claim that when $\phi \in \mathcal{A}_N$, the grand maximal function $\mathscr{M}_N(f)$ satisfies

$$|\mathscr{M}_N(f)(x)| \leqslant C\int_{\mathbb{R}^n} |f(x)|/|x-y|^n dy. \tag{2}$$

In fact,

$$|\mathscr{M}_N(f)(x)| = \sup_t \sup_{|x-y|<t} |\phi_t * f(y)|$$

$$\leqslant \sup_t \sup_{|x-y|<t} \int_{\mathbb{R}^n} |\phi_t(y-x)| |f(y)| dy$$

$$\lesssim \sup_t \sup_{|x-y|<t} \int_{\mathbb{R}^n} \left| t^{-n} \left( \frac{t}{y-x} \right)^n \right| |f(y)|\, dy$$

$$\lesssim \int_{\mathbb{R}^n} \frac{|f(x)|}{|y-x|^n} dy.$$

It follows from (2) and Lemma 2.1 that

$$\|f\|_{H\dot{K}_{\vec{p}}^{\alpha,q}(\mathbb{R}^n)} = \|\mathscr{M}_N(f)\|_{\dot{K}_{\vec{p}}^{\alpha,q}(\mathbb{R}^n)} \leqslant C \, \|f\|_{\dot{K}_{\vec{p}}^{\alpha,q}(\mathbb{R}^n)}.$$

As a consequence, $\dot{K}_{\vec{p}}^{\alpha,q}(\mathbb{R}^n) \subset H\dot{K}_{\vec{p}}^{\alpha,q}(\mathbb{R}^n) \bigcap L_{\text{loc}}^{\vec{p}}(\mathbb{R}^n \backslash \{0\})$.

In addition, we can easily obtain $\dot{K}_{\vec{p}}^{\alpha,q}(\mathbb{R}^n) \subset L_{\text{loc}}^{\vec{p}}(\mathbb{R}^n)$ from the definition of mixed Herz spaces.
□

**Remark 2.3.** *In view of Proposition 2.3, we will only consider the characterizations of the spaces $H\dot{K}_{\vec{p}}^{\alpha,q}(\mathbb{R}^n)$ and $HK_{\vec{p}}^{\alpha,q}(\mathbb{R}^n)$ in terms of central atom and molecular decomposition for $\alpha \geqslant n(1 - 1/n \sum_{i=1}^m 1/p_i)$ in the subsequent sections.*

**Proposition 2.4.** *Let $\alpha \in \mathbb{R}, 0 < q < \infty, 1 < \vec{p} < \infty$. Then*

$$HK_{\vec{p}}^{\alpha,q}(\mathbb{R}^n) = L^{\vec{p}}(\mathbb{R}^n) \bigcap H\dot{K}_{\vec{p}}^{\alpha,q}(\mathbb{R}^n).$$

*Proof.* To prove this proposition we need the following assertion. For $0 < q, \vec{p} \leqslant \infty$, we have

$$K_{\vec{p}}^{\alpha,q}(\mathbb{R}^n) = \dot{K}_{\vec{p}}^{\alpha,q}(\mathbb{R}^n) \bigcap L^{\vec{p}}(\mathbb{R}^n), \tag{3}$$

and for all $f \in \dot{K}_{\vec{p}}^{\alpha,q}(\mathbb{R}^n) \bigcap L^{\vec{p}}(\mathbb{R}^n)$, there holds

$$\|f\|_{K_{\vec{p}}^{\alpha,q}(\mathbb{R}^n)} = \|f\|_{\dot{K}_{\vec{p}}^{\alpha,q}(\mathbb{R}^n)} + \|f\|_{L^{\vec{p}}(\mathbb{R}^n)}.$$

Note that

$$\|f\|_{K_{\vec{p}}^{\alpha,q}(\mathbb{R}^n)}^q = \sum_{k=0}^{\infty} 2^{k\alpha q} \|f\widetilde{\chi}_k\|_{L^{\vec{p}}(\mathbb{R}^n)}^q$$

$$= \|f\widetilde{\chi}_0\|_{L^{\vec{p}}(\mathbb{R}^n)}^q + \sum_{k=1}^{\infty} 2^{k\alpha q} \|f\chi_k\|_{L^{\vec{p}}(\mathbb{R}^n)}^q$$

$$\leqslant \|f\|_{L^{\vec{p}}(\mathbb{R}^n)}^q + \|f\|_{\dot{K}_{\vec{p}}^{\alpha,q}(\mathbb{R}^n)}^q - \sum_{k=-\infty}^{0} 2^{k\alpha q} \|f\chi_k\|_{L^{\vec{p}}(\mathbb{R}^n)}^q.$$

We claim that the following estimation is correct:

$$\sum_{k=-\infty}^{0} 2^{k\alpha q} \|f\chi_k\|_{L^{\vec{p}}(\mathbb{R}^n)}^q \leqslant C \|f\|_{L^{\vec{p}}(|x|\leqslant 1)}^q.$$

When $q \leqslant p_n$, by using Hölder's inequality and

$$\left( \sum_{k=1}^{\infty} |a_k| \right)^r \leqslant \sum_{k=1}^{\infty} |a_k|^r \qquad (0 < r \leqslant 1),$$

we obtain

$$\sum_{k=-\infty}^{0} 2^{k\alpha q}\|f\chi_k\|_{L^{\vec{p}}(\mathbb{R}^n)}^{q} \leqslant \left(\sum_{k=-\infty}^{0} 2^{k\alpha q(\frac{q}{p_n})'}\right)^{\frac{1}{(\frac{p_n}{q})'}} \left(\sum_{k=-\infty}^{0} \|f\chi_k\|_{L^{\vec{p}}(\mathbb{R}^n)}^{q\frac{p_n}{q}}\right)^{\frac{q}{p_n}} \leqslant C\|f\|_{L^{\vec{p}}(|x|\leqslant 1)}^{q}.$$

When $p_n < q$,

$$\sum_{k=-\infty}^{0} 2^{k\alpha q}\|f\chi_k\|_{L^{\vec{p}}(\mathbb{R}^n)}^{q} = \sum_{k=-\infty}^{0} \left(2^{k\alpha p_n}\|f\chi_k\|_{L^{\vec{p}}(\mathbb{R}^n)}^{p_n}\right)^{\frac{q}{p_n}} \leqslant \left(C\|f\|_{L^{\vec{p}}(|x|\leqslant 1)}^{p_n}\right)^{\frac{q}{p_n}} \leqslant C\|f\|_{L^{\vec{p}}(|x|\leqslant 1)}^{q},$$

we also obtain $\|f\|_{K^{\alpha,q}_{\vec{p}}(\mathbb{R}^n)} \leqslant \|f\|_{\dot{K}^{\alpha,q}_{\vec{p}}(\mathbb{R}^n)} + \|f\|_{L^{\vec{p}}(\mathbb{R}^n)}$.

To finish the proof, we need to check $\|f\|_{\dot{K}^{\alpha,q}_{\vec{p}}(\mathbb{R}^n)} \leqslant C\|f\|_{K^{\alpha,q}_{\vec{p}}(\mathbb{R}^n)}$ and $\|f\|_{L^{\vec{q}}(\mathbb{R}^n)} \leqslant C\|f\|_{K^{\alpha,p}_{\vec{q}}(\mathbb{R}^n)}$.

By a direct computation, we have

$$\|f\|_{\dot{K}^{\alpha,q}_{\vec{p}}(\mathbb{R}^n)}^{q} \leqslant \sum_{k=-\infty}^{0} 2^{k\alpha q}\|f\chi_k\|_{L^{\vec{p}}(\mathbb{R}^n)}^{q} + \sum_{k=1}^{\infty} 2^{k\alpha q}\|f\chi_k\|_{L^{\vec{p}}(\mathbb{R}^n)}^{q}$$

$$\leqslant C\|f\|_{L^{\vec{p}}(|x|\leqslant 1)}^{q} + \sum_{k=0}^{\infty} 2^{k\alpha q}\|f\chi_k\|_{L^{\vec{p}}(\mathbb{R}^n)}^{q}$$

$$\leqslant C\|f\|_{L^{\vec{p}}(|x|\leqslant 1)}^{q} + \|f\|_{K^{\alpha,q}_{\vec{p}}(\mathbb{R}^n)}.$$

The result $\|f\|_{L^{\vec{p}}(|x|\leqslant 1)} \leqslant C\|f\|_{K^{\alpha,q}_{\vec{p}}(\mathbb{R}^n)}$ is evident.

Furthermore, to prove the other estimate, it suffices to show that

$$\|f\|_{L^{\vec{p}}(|x|>1)} \leqslant C\|f\|_{K^{\alpha,q}_{\vec{p}}(\mathbb{R}^n)}.$$

When $0 < p_n \leqslant q$, as a consequence of Hölder's inequality, there holds

$$\|f\|_{L^{\vec{p}}(|x|>1)}^{p_n} \leqslant \left(\sum_{k=-\infty}^{0} 2^{-k\alpha p_n(\frac{q}{p_n})'}\right)^{\frac{1}{(\frac{q}{p_n})'}} \left(\sum_{k=-\infty}^{0} 2^{k\alpha q}\|f\chi_k\|_{L^{\vec{p}}(\mathbb{R}^n)}^{p_n\frac{q}{p_n}}\right)^{\frac{p_n}{q}} \leqslant C\|f\|_{K^{\alpha,q}_{\vec{p}}(\mathbb{R}^n)}^{p_n}.$$

When $p_n > q$, we also have

$$\|f\|_{L^{\vec{p}}(|x|>1)}^{p_n} \leqslant \sum_{k=1}^{\infty} \|f\chi_k\|_{L^{\vec{p}}(\mathbb{R}^n)}^{p_n} \leqslant \left(\sum_{k=1}^{\infty} \|f\chi_k\|_{L^{\vec{p}}(\mathbb{R}^n)}^{q}\right)^{\frac{p_n}{q}} \leqslant C\|f\|_{K^{\alpha,q}_{\vec{p}}(\mathbb{R}^n)}^{p_n}.$$

The above estimates indicate that $K^{\alpha,q}_{\vec{p}}(\mathbb{R}^n) = \dot{K}^{\alpha,q}_{\vec{p}}(\mathbb{R}^n) \bigcap L^{\vec{p}}(\mathbb{R}^n)$ holds.

By using (3) and the boundedness of $\mathscr{M}_N$ on mixed spaces $L^{\vec{p}}(\mathbb{R}^n)$ from Proposition 2.2, we can immediately show this proposition. □

## 3 Atom characterization of mixed Herz-Hardy spaces

We begin with the definition of central atoms.

**Definition 3.1.** *Let $1 < \vec{p} < \infty$, $n(1 - 1/n\sum_{i=1}^{m} 1/p_i) \leqslant \alpha < \infty$, and non-negative integer $s \geqslant [\alpha - n(1 - 1/n\sum_{i=1}^{m} 1/p_i)]$.*

(i) A function $a$ on $\mathbb{R}^n$ is said to be a central $(\alpha, \vec{p})$-atom, if it satisfies
(a) $\operatorname{supp} a \subset B(0, r) = \{x \in \mathbb{R}^n : |x| < r\}$.
(b) $\|a\|_{L^{\vec{p}}(\mathbb{R}^n)} \leqslant |B(0, r)|^{-\alpha/n}$.
(c) $\int_{\mathbb{R}^n} a(x) x^\beta dx = 0, |\beta| \leqslant s$.
(ii) A function $a$ on $\mathbb{R}^n$ is said to be a central $(\alpha, \vec{p})$-atom of restricted type, if it satisfies the conditions (b), (c) above and
(a)' $\operatorname{supp} a \subset B(0, r), r \geqslant 1$.

**Remark 3.1.** If $\vec{p} = p$ is a constant, we recover the classical $(\alpha, p)$-atom since $1 - 1/n \sum_{i=1}^{m} 1/p_i = 1 - 1/p$.

**Theorem 3.1.** Let $1 < \vec{p} < \infty$, $0 < q < \infty$ and $n(1 - 1/n \sum_{i=1}^{m} 1/p_i) \leqslant \alpha < \infty$. Then
(i) $f \in H\dot{K}_{\vec{p}}^{\alpha, q}(\mathbb{R}^n)$ if and only if

$$f = \sum_{k=-\infty}^{\infty} \lambda_k a_k \quad \text{in the sense of } \mathcal{S}'(\mathbb{R}^n),$$

where each $a_k$ is a central $(\alpha, \vec{p})$-atom with support contained in $B_k$ and $\sum_{k=0}^{\infty} |\lambda_k|^q < \infty$. Moreover,

$$\|f\|_{H\dot{K}_{\vec{p}}^{\alpha, q}(\mathbb{R}^n)} \approx \inf \left( \sum_{k=0}^{\infty} |\lambda_k|^q \right)^{1/q},$$

where the infimum is taken over all above decompositions of $f$.
(ii) $f \in HK_{\vec{p}}^{\alpha, q}(\mathbb{R}^n)$ if and only if

$$f = \sum_{k=0}^{\infty} \lambda_k a_k \quad \text{in the sense of } \mathcal{S}'(\mathbb{R}^n),$$

where each $a_k$ is a central $(\alpha, \vec{p})$-atom of restricted type with support contained in $B_k$ and $\sum_{k=0}^{\infty} |\lambda_k|^q < \infty$. Moreover,

$$\|f\|_{HK_{\vec{p}}^{\alpha, q}(\mathbb{R}^n)} \approx \inf \left( \sum_{k=0}^{\infty} |\lambda_k|^q \right)^{1/q},$$

where the infimum is taken over all above decompositions of $f$.

*Proof.* We just need to show (i), and (ii) can be proved in the similar way.
Necessity: choose $\phi \in C_0^\infty(\mathbb{R}^n)$ with $\phi \geqslant 0$, $\int_{\mathbb{R}^n} \phi(x) dx = 1$, $\operatorname{supp} \phi \subset \{x : |x| \leqslant 1\}$. For $j \in \mathbb{Z}_+$, let

$$\phi_{(j)}(x) = 2^{jn} \phi(2^j x).$$

For each $f \in \mathcal{S}'(\mathbb{R}^n)$, write
$$f^{(j)}(x) = f * \phi_{(j)}(x).$$

It is clearly that $f^{(j)} \in C^\infty(\mathbb{R}^n)$ and $\lim_{j \to \infty} f^{(j)} = f$ in the sense of distribution.

Let $\psi$ be a radial smooth function such that $\operatorname{supp} \psi \subset \{x : 1/2 - \varepsilon \leqslant |x| \leqslant 1 + \varepsilon\}$ with $0 < \varepsilon < 1/4$, $\psi(x) = 1$ for $1/2 \leqslant |x| \leqslant 1$. Set $\psi_k(x) = \psi(2^{-k} x)$ for $k \in \mathbb{Z}$ and

$$\tilde{A}_{k, \varepsilon} = \left\{ x : 2^{k-1} - 2^k \varepsilon \leqslant |x| \leqslant 2^k + 2^k \varepsilon \right\}.$$

From the conditions of $\psi$, we know that $\operatorname{supp} \psi_k \subset \tilde{A}_{k,\varepsilon}$ and $\psi_k(x) = 1$ for $x \in A_k = \{x : 2^{k-1} < |x| \leqslant 2^k\}$.

Obviously, $1 \leqslant \sum_{k=-\infty}^{\infty} \psi_k(x) \leqslant 2, |x| > 0$. Write

$$\Phi_k(x) = \begin{cases} \psi_k(x) / \sum_{l=-\infty}^{\infty} \psi_l(x), & x \neq 0, \\ 0, & x = 0, \end{cases}$$

then $\sum_{k=-\infty}^{\infty} \Phi_k(x) = 1$ for $x \neq 0$. For some $m \in \mathbb{N}$, we denote the class of all the real polynomials with the degree less than $m$ by $\mathcal{P}_m(\mathbb{R}^n)$. Let $P_k^{(j)}(x) = P_{\tilde{A}_{k,\varepsilon}}\left(f^{(j)}\Phi_k\right)(x)\chi_{\tilde{A}_{k,\varepsilon}}(x) \in \mathcal{P}_m(\mathbb{R}^n)$ be the unique polynomial satisfying

$$\int_{\tilde{A}_{k,\varepsilon}} \left(f^{(j)}(x)\Phi_k(x) - P_k^{(j)}(x)\right) x^\beta dx = 0, \quad |\beta| \leqslant m = \left[\alpha - n\left(1 - \frac{1}{n}\sum_{i=1}^{n}\frac{1}{p_i}\right)\right].$$

Denote by

$$f^{(j)}(x) = \sum_{k=-\infty}^{\infty} \left(f^{(j)}(x)\Phi_k(x) - P_k^{(j)}(x)\right) + \sum_{k=-\infty}^{\infty} P_k^{(j)}(x) = \sum_{I}^{(j)} + \sum_{II}^{(j)}.$$

For the term $\sum_{I}^{(j)}$, let $g_k^{(j)}(x) = f^{(j)}(x)\Phi_k(x) - P_k^{(j)}(x)$.

Now we deal with $\left\|g_k^{(j)}\right\|_{L^{\vec{p}}(\mathbb{R}^n)}$. To do this, we need to estimate $\left|P_k^{(j)(x)}\right|$.

Let $\{\phi_d^k : |d| \leqslant m\}$ be the orthogonal polynomials restricted to $\tilde{A}_{k,\varepsilon}$ with respect to the weight $1/|\tilde{A}_{k,\varepsilon}|$, which are obtained from $\{x^\beta : |\beta| \leqslant m\}$ by Gram-Schmidt's method, that is

$$\left\langle \phi_\nu^k, \phi_\mu^k \right\rangle = \frac{1}{|\tilde{A}_{k,\varepsilon}|} \int_{\tilde{A}_{k,\varepsilon}} \phi_\nu^k(x)\phi_\mu^k(x)dx = \delta_{\nu\mu}.$$

We know that $P_k^{(j)}(x) = \sum_{|d|\leqslant m} \langle f^{(j)}\Phi_k, \phi_d^k \rangle \phi_d^k(x)$ for $x \in \tilde{A}_{k,\varepsilon}$. On the other hand, from dilation, we can gain that

$$\frac{1}{|\tilde{A}_{1,\varepsilon}|} \int_{\tilde{A}_{1,\varepsilon}} \phi_\nu^k\left(2^{k-1}y\right) \phi_\mu^k\left(2^{k-1}y\right) dy = \delta_{\nu\mu}.$$

Moreover, there holds $\phi_\nu^k(x) = \phi_\nu^1(2^{1-k}x)$ a.e. for $x \in \tilde{A}_{k,\varepsilon}$. Thus $|\phi_\nu^k(x)| \leqslant C$, by Hölder's inequality on mixed norm spaces, we have

$$\left|P_k^{(j)}(x)\right| \leqslant \frac{C}{|\tilde{A}_{k,\varepsilon}|} \int_{\tilde{A}_{k,\varepsilon}} \left|f^{(j)}(x)\Phi_k(x)\right| dx$$

$$\leqslant \frac{C}{|\tilde{A}_{k,\varepsilon}|} \left\|f^{(j)}\Phi_k\right\|_{L^{\vec{p}}(\mathbb{R}^n)} \left\|\chi_{\tilde{A}_{k,\varepsilon}}\right\|_{L^{\vec{p}'}(\mathbb{R}^n)}.$$

Therefore, we have

$$\left\|g_k^{(j)}\right\|_{L^{\vec{p}}(\mathbb{R}^n)} \leqslant \left\|f^{(j)}\Phi_k\right\|_{L^{\vec{p}}(\mathbb{R}^n)} + \left\|P_k^{(j)}\right\|_{L^{\vec{p}}(\mathbb{R}^n)}$$

$$\leqslant \left\|f^{(j)}\Phi_k\right\|_{L^{\vec{p}}(\mathbb{R}^n)} + \frac{C}{\left|\tilde{A}_{k,\varepsilon}\right|}\left\|f^{(j)}\Phi_k\right\|_{L^{\vec{p}}(\mathbb{R}^n)}\left\|\chi_{\tilde{A}_{k,\varepsilon}}\right\|_{L^{\vec{p'}}(\mathbb{R}^n)}\left\|\chi_{\tilde{A}_{k,\varepsilon}}\right\|_{L^{\vec{p}}(\mathbb{R}^n)}$$

$$\leqslant \left\|f^{(j)}\Phi_k\right\|_{L^{\vec{p}}(\mathbb{R}^n)} + C\left\|f^{(j)}\Phi_k\right\|_{L^{\vec{p}}(\mathbb{R}^n)}$$

$$\leqslant C\left\|\left(f*\phi_{(j)}\right)\Phi_k\right\|_{L^{\vec{p}}(\mathbb{R}^n)}$$

$$\leqslant C\sum_{l=k-1}^{k+1}\left\|(\mathscr{M}_N f)\chi_l\right\|_{L^{\vec{p}}(\mathbb{R}^n)}.$$

Now write

$$\sum_{k=-\infty}^{\infty}\left(f^{(j)}(x)\Phi_k(x) - P_k^{(j)}(x)\right) = \sum_{k=-\infty}^{\infty} |B_{k+1}|^{\frac{\alpha}{n}} \sum_{l=k-1}^{k+1} \|(\mathscr{M}_N f)\chi_l\|_{L^{\vec{p}}(\mathbb{R}^n)}$$

$$\times \frac{f^{(j)}(x)\Phi_k(x) - P_k^{(j)}(x)}{|B_{k+1}|^{\frac{\alpha}{n}}\sum_{l=k-1}^{k+1}\|(\mathscr{M}_N f)\chi_l\|_{L^{\vec{p}}(\mathbb{R}^n)}} = \sum_{k=-\infty}^{\infty} \lambda_k a_k^{(j)}.$$

Then $\left\|a_k^{(j)}\right\|_{L^{\vec{p}}(\mathbb{R}^n)} \leqslant |B_{k+1}|^{-\alpha/n}$ and each $a_k^{(j)}$ is a central $(\alpha, \vec{p})$-atom with support contained in $B_{k+1}$. Furthermore,

$$\sum_{k=-\infty}^{\infty} |\lambda_k|^q \leqslant C \sum_{k=-\infty}^{\infty} |B_{k+1}|^{\frac{\alpha q}{n}} \sum_{l=k-1}^{k+1} \|(\mathscr{M}_N(f)\chi_l)\|_{L^{\vec{p}}(\mathbb{R}^n)}^q \leqslant C\|\mathscr{M}_N f\|_{\dot{K}_{\vec{p}}^{\alpha,q}(\mathbb{R}^n)}^q \leqslant C\|f\|_{H\dot{K}_{\vec{p}}^{\alpha,q}(\mathbb{R}^n)}^q.$$

It remains to estimate $\sum_{II}^{(j)}$. Let $\{\psi_d^k : |d| \leqslant m\}$ be the dual basis of $\{x^\beta : |\beta| \leqslant m\}$ with respect to the weight $1/\left|\tilde{A}_{k,\varepsilon}\right|$ on $\tilde{A}_{k,\varepsilon}$, that is

$$\left\langle \psi_d^k, x^\beta \right\rangle = \frac{1}{\left|\tilde{A}_{k,\varepsilon}\right|} \int_{\tilde{A}_{k,\varepsilon}} x^\beta \psi_d^k(x) dx = \delta_{\beta d}.$$

Let

$$h_{k,d}^{(j)}(x) = \sum_{l=-\infty}^{k} \left(\frac{\psi_d^k(x)\chi_{\tilde{A}_{k,\varepsilon}}(x)}{\left|\tilde{A}_{k,\varepsilon}\right|} - \frac{\psi_d^{k+1}(x)\chi_{\tilde{A}_{k+1,\varepsilon}}(x)}{\left|\tilde{A}_{k+1,\varepsilon}\right|}\right) \int_{\mathbb{R}^n} f^{(j)}(x)\Phi_l(x)x^d dx.$$

We can write

$$\sum_{II}^{(j)} = \sum_{|d|\leqslant m}\sum_{k=-\infty}^{\infty}\left(\int_{\mathbb{R}^n} f^{(j)}\Phi_k x^d dx\right)\frac{\psi_d^k(x)\chi_{\tilde{A}_{k,\varepsilon}}(x)}{\left|\tilde{A}_{k,\varepsilon}\right|}$$

$$= \sum_{|d|\leqslant m}\sum_{k=-\infty}^{\infty}\left(\sum_{l=-\infty}^{k}\int_{\mathbb{R}^n} f^{(j)}(x)\Phi_l(x)x^d dx\right)$$

$$\times \left(\frac{\psi_d^k(x)\chi_{\tilde{A}_{k,\varepsilon}}(x)}{\left|\tilde{A}_{k,\varepsilon}\right|} - \frac{\psi_d^{k+1}(x)\chi_{\tilde{A}_{k+1,\varepsilon}}(x)}{\left|\tilde{A}_{k+1,\varepsilon}\right|}\right)$$

$$= \sum_{|d|\leqslant m} \sum_{k=-\infty}^{\infty} \frac{\alpha_{k,d} h_{k,d}^{(j)}(x)}{\alpha_{k,d}} = \sum_{|d|\leqslant m} \sum_{k=-\infty}^{\infty} \alpha_{k,d} a_{k,d}^{(j)}(x),$$

where

$$\alpha_{k,d} = C \sum_{l=k-1}^{k+1} \|(\mathscr{M}_N f)\chi_l\|_{L^{\vec{p}}(\mathbb{R}^n)} |B_{k+2}|^{\frac{\alpha}{n}}.$$

Note that

$$\int_{\mathbb{R}^n} \sum_{l=-\infty}^{k} \left|\Phi_l(x)x^d\right| dx = \sum_{l=-\infty}^{k} \int_{\tilde{A}_{k,\varepsilon}} \left|\Phi_l(x)x^d\right| dx \leqslant C 2^{k(n+|d|)}.$$

Furthermore,

$$\left|\int_{\mathbb{R}^n} f^{(j)}(y) \sum_{l=-\infty}^{k} \Phi_l(y) y^d dy\right| \leqslant C 2^{k(n+|d|)} \mathscr{M}_N f(x), \quad x \in B_{k+2}. \tag{4}$$

Inequality (4), together with the inequality that

$$\left(\frac{\psi_d^k(x)\chi_{\tilde{A}_{k,\varepsilon}}(x)}{\left|\tilde{A}_{k,\varepsilon}\right|} - \frac{\psi_d^{k+1}(x)\chi_{\tilde{A}_{k+1,\varepsilon}}(x)}{\left|\tilde{A}_{k+1,\varepsilon}\right|}\right) \leqslant C 2^{-k(n+|d|)} \sum_{l=k-1}^{k+1} \chi_l(x),$$

shows that

$$\left\|h_{k,d}^{(j)}\right\|_{L^{\vec{p}}(\mathbb{R}^n)} \leqslant C \sum_{l=k-1}^{k+1} \|(\mathscr{M}_N f)\chi_l\|_{L^{\vec{p}}(\mathbb{R}^n)},$$

and

$$\left\|a_{k,d}^{(j)}\right\|_{L^{\vec{p}}(\mathbb{R}^n)} \leqslant \left\|\frac{h_{k,d}^{(j)}}{C \sum_{l=k-1}^{k+1} \|(\mathscr{M}_N f)\chi_l\|_{L^{\vec{p}}(\mathbb{R}^n)} |B_{k+2}|^{-\frac{\alpha}{n}}}\right\|_{L^{\vec{p}}(\mathbb{R}^n)} \leqslant C |B_{k+2}|^{-\frac{\alpha}{n}}.$$

It has been check that $a_{k,d}^{(j)}$ is a central $(\alpha, \vec{p})$-atom with support contained in $\tilde{A}_{k,\varepsilon} \cup \tilde{A}_{k+1,\varepsilon} \subset B_{k+2}$, and $\alpha_{k,d} = C \sum_{l=k-1}^{k+2} \|(Gf)\chi_l\|_{L^{\vec{p}}(\mathbb{R}^n)} |B_{k+2}|^{\alpha/n}$. Moreover, we also can get

$$\sum_{k,d} |\alpha_{k,d}|^q \leqslant C \sum_{k=-\infty}^{\infty} |B_{k+2}|^{\frac{\alpha q}{n}} \left(\sum_{l=k-1}^{k+1} \|(\mathscr{M}_N f)\chi_l\|_{L^{\vec{p}}(\mathbb{R}^n)}\right)^q \leqslant C \|\mathscr{M}_N f\|_{K_{\vec{p}}^{\alpha,q}(\mathbb{R}^n)}^q < \infty.$$

Thus we obtain that

$$f^{(j)}(x) = \sum_{d=-\infty}^{\infty} \lambda_d a_d^{(j)}(x),$$

where each $a_d^{(j)}$ is a central $(\alpha, \vec{p})$-atom with support contained in $\tilde{A}_{d,\varepsilon} \cup \tilde{A}_{d+1,\varepsilon} \subset B_{d+2}$, and

$$\left(\sum_{d=-\infty}^{\infty} |\lambda_d|^q\right)^{1/q} \leqslant C \|\mathscr{M}_N f\|_{\dot{K}_{\vec{p}}^{\alpha,p}(\mathbb{R}^n)}^q < \infty.$$

Since $\sup_{j\in\mathbb{Z}_+}\left\|a_0^{(j)}\right\|_{L^{\vec{p}}(\mathbb{R}^n)}\leqslant |B_2|^{-\alpha/n}$, by the Banach-Alaoglu theorem we can obtain a subsequence $\left\{a_0^{(j_{n_0})}\right\}$ of $\left\{a_0^{(j)}\right\}$ converging in the weak* topology of $L^{\vec{p}}(\mathbb{R}^n)$ to some $a_0 \in L^{\vec{p}}(\mathbb{R}^n)$ and $a_0$ is a central $(\alpha,\vec{p})$-atom supported on $B_2$. Repeating the above procedure for each $d\in\mathbb{Z}$, we can find a subsequence $\left\{a_d^{(j_{n_d})}\right\}$ of $\left\{a_d^{(j)}\right\}$ converging weak $*$ in $L^{\vec{p}}(\mathbb{R}^n)$ to some $a_d \in L^{\vec{p}}(\mathbb{R}^n)$ which is a central $(\alpha,\vec{p})$-atom supported on $B_{d+2}$. By the usual diagonal method we obtain a subsequence $\{j_\nu\}$ of $\mathbb{Z}_+$ such that for each $d\in\mathbb{Z}, \lim_{\nu\to\infty} a_d^{(j_\nu)} = a_d$ in the weak* topology of $L^{\vec{p}}(\mathbb{R}^n)$ and in $\mathcal{S}'(\mathbb{R}^n)$.

Now our proof is reduced to prove that

$$f = \sum_{d=-\infty}^{\infty} \lambda_d a_d, \text{ in the sense of } \mathcal{S}'(\mathbb{R}^n).$$

For each $\varphi \in \mathcal{S}(\mathbb{R}^n)$, noting that $\operatorname{supp} a_d^{(j_\nu)} \subset \left(\tilde{A}_{d,\varepsilon} \cup \tilde{A}_{d+1,\varepsilon}\right) \subset (A_{d-1}\cup A_d \cup A_{d+1} \cup A_{d+2})$, we have

$$\langle f, \varphi\rangle = \lim_{\nu\to\infty}\sum_{d=-\infty}^{\infty} \lambda_d \int_{\mathbb{R}^n} a_d^{(j_\nu)}(x)\varphi(x)dx.$$

Recall that $m = [\alpha - \sum_{i=1}^n 1/p_i']$. If $d \leqslant 0$, by using Hölder's inequality on mixed-norm spaces, then

$$\left|\int_{\mathbb{R}^n} a_d^{(j_\nu)}(x)\varphi(x)dx\right| = \left|\int_{\mathbb{R}^n} a_d^{(j_\nu)}(x)\left(\varphi(x) - \sum_{|\beta|\leqslant m}\frac{D^\beta\varphi(0)}{\beta!}x^\beta\right)dx\right|$$

$$\leqslant C\int_{\mathbb{R}^n}\left|a_d^{(j_\nu)}(x)\right|\cdot|x|^{m+1}dx$$

$$\leqslant C 2^{d(m+1)}\left\|a_d^{(j_\nu)}\right\|_{L^{\vec{p}}(\mathbb{R}^n)}\left\|\chi_{B_{d+2}}\right\|_{L^{\vec{p}'}(\mathbb{R}^n)}$$

$$\leqslant C 2^{d\left(m+1-\alpha+\sum_{i=1}^n\frac{1}{p_i'}\right)}.$$

If $d > 0$, let $k_0 \in \mathbb{Z}_+$ such that $k_0 \geqslant \sum_{i=1}^n 1/p_i' - \alpha$. Again by using Hölder's inequality on mixed-norm spaces, it yields

$$\left|\int_{\mathbb{R}^n} a_d^{(j_\nu)}(x)\varphi(x)dx\right| \leqslant C\int_{\mathbb{R}^n}\left|a_d^{(j_\nu)}(x)\right||x|^{-k_0}dx$$

$$\leqslant C 2^{-dk_0}\left\|a_d^{(j_\nu)}\right\|_{L^{\vec{p}}(\mathbb{R}^n)}\left\|\chi_{B_{d+2}}\right\|_{L^{\vec{p}'}(\mathbb{R}^n)}$$

$$\leqslant C 2^{-d\left(k_0+\alpha-\sum_{i=1}^n\frac{1}{p_i'}\right)}.$$

Let

$$\mu_d = \begin{cases} |\lambda_d|\,2^{d\left(m+1-\alpha+\sum_{i=1}^n\frac{1}{p_i'}\right)}, & d\leqslant 0, \\ |\lambda_d|\,2^{-d\left(k_0+\alpha-\sum_{i=1}^n\frac{1}{p_i'}\right)}, & d > 0. \end{cases}$$

Then
$$\sum_{d=-\infty}^{\infty} |\mu_d| \leqslant C \left( \sum_{d=-\infty}^{\infty} |\lambda_d|^q \right)^{\frac{1}{q}} \leqslant C \|\mathscr{M}_N f\|_{\dot{K}_{\vec{p}}^{\alpha,q}(\mathbb{R}^n)} < \infty$$
and
$$|\lambda_d| \left| \int_{\mathbb{R}^n} a_d^{(j_\nu)}(x) \varphi(x) dx \right| \leqslant C |\mu_d|,$$
which imply that
$$\langle f, \varphi \rangle = \sum_{d=-\infty}^{\infty} \lim_{\nu \to \infty} \lambda_d \int_{\mathbb{R}^n} a_d^{(j_\nu)}(x) \varphi(x) dx = \sum_{d=-\infty}^{\infty} \lambda_d \int_{\mathbb{R}^n} a_d(x) \varphi(x) dx. \tag{5}$$

This mean that $f = \sum_{d=-\infty}^{\infty} \lambda_d a_d$ in the sense of $\mathcal{S}'(\mathbb{R}^n)$.

Sufficiency: we will prove the conclusion for two different cases: $0 < q \leqslant 1$ and $1 < q < \infty$.

If $0 < q \leqslant 1$, it suffices to show that for each central $(\alpha, \vec{p})$-atom $a$,
$$\|\mathscr{M}_N a\|_{\dot{K}_{\vec{p}}^{\alpha,q}(\mathbb{R}^n)} \leqslant C$$
with the constant $C > 0$ independent of $a$.

For a fixed central $(\alpha, \vec{p})$-atom $a$, with $\operatorname{supp} a(x) \subset B(0, 2^{k_0})$ for some $k_0 \in \mathbb{Z}$. Write
$$\|\mathscr{M}_N a\|_{\dot{K}_{\vec{p}}^{\alpha,q}(\mathbb{R}^n)}^q = \sum_{k=-\infty}^{k_0+3} 2^{k\alpha q} \|(\mathscr{M}_N a) \chi_k\|_{L^{\vec{p}}(\mathbb{R}^n)}^q$$
$$+ \sum_{k=k_0+4}^{\infty} 2^{k\alpha q} \|(\mathscr{M}_N a) \chi_k\|_{L^{\vec{p}}(\mathbb{R}^n)}^q$$
$$= I + II.$$

By the $L^{\vec{p}}(\mathbb{R}^n)$ boundedness of the grand maximal operator $\mathscr{M}_N$ from Proposition 2.2, we have
$$I \leqslant \sum_{k=-\infty}^{k_0+3} 2^{k\alpha q} \|\mathscr{M}_N a\|_{L^{\vec{p}}(\mathbb{R}^n)}^q \leqslant C \|a\|_{L^{\vec{p}}(\mathbb{R}^n)}^q \sum_{k=-\infty}^{k_0+3} 2^{k\alpha q} \leqslant C.$$

The next step is to consider $II$. We need a pointwise estimate for $\mathscr{M}_N a(x)$ on $A_k$.

Let $\phi \in \mathcal{A}_N$, $m \in \mathbb{N}$ such that $\alpha - \sum_{i=1}^{n} \frac{1}{p_i'} < m+1$. Denote by $P_m$ the $m$-th order Taylor series expansion. If $|x - y| < t$, then from the vanishing moment condition of $a$ we can get
$$|a * \phi_t(y)| = t^{-n} \left| \int_{\mathbb{R}^n} a(z) \left( \phi\left(\frac{y-z}{t}\right) - P_m\left(\frac{y}{t}\right) \right) dz \right|$$
$$\leqslant C t^{-n} \int_{\mathbb{R}^n} |a(z)| \left|\frac{z}{t}\right|^{m+1} \left(1 + \frac{|y - \theta z|}{t}\right)^{-(n+m+1)} dz$$
$$\leqslant C \int_{\mathbb{R}^n} |a(z)| |z|^{m+1} (t + |y - \theta z|)^{-(n+m+1)} dz,$$
where $0 < \theta < 1$.

Since $x \in A_k$ for $k \geqslant k_0+4$, we have $|x| \geqslant 2 \cdot 2^{k_0+1}$. From $|x-y| < t$ and $|z| < 2^{k_0+1}$, we have

$$t + |y - \theta z| \geqslant |x-y| + |y - \theta z| \geqslant |x| - |z| \geqslant |x|/2.$$

Thus,

$$\begin{aligned}|a * \phi_t(y)| &\leqslant C \int_{\mathbb{R}^n} |a(z)||z|^{m+1}(|x-y| + |y - \theta z|)^{-(n+m+1)} dz \\ &\leqslant C 2^{k_0(m+1)} |x|^{-(n+m+1)} \int_{\mathbb{R}^n} |a(z)| dz \\ &\leqslant C 2^{k_0(m+1)} |x|^{-(n+m+1)} |B_{k_0}|^{-\frac{\alpha}{n}} \left\| \chi_{B_{k_0}} \right\|_{L^{\vec{p}'}(\mathbb{R}^n)}.\end{aligned}$$

Therefore, we have

$$\mathscr{M}_N a(x) \leqslant C 2^{k_0(m+1) - k(n+m+1)} |B_{k_0}|^{-\frac{\alpha}{n}} \left\| \chi_{B_{k_0}} \right\|_{L^{\vec{p}'}(\mathbb{R}^n)}, x \in A_k \quad k \geqslant k_0 + 4.$$

As a consequence,

$$\begin{aligned}II &= \sum_{k=k_0+4}^{\infty} 2^{k\alpha q} \left\| (\mathscr{M}_N a) \chi_k \right\|_{L^{\vec{p}}(\mathbb{R}^n)}^q \\ &\leqslant C \sum_{k=k_0+4}^{\infty} 2^{k\alpha q} 2^{k_0(m+1) - k(n+m+1)} 2^{k_0 \alpha q} \left\| \chi_{B_{k_0}} \right\|_{L^{\vec{p}'}(\mathbb{R}^n)}^q \| \chi_{B_k} \|_{L^{\vec{p}}(\mathbb{R}^n)}^q \\ &\leqslant C \sum_{k=k_0+4}^{\infty} 2^{k_0(m+1-\alpha) - k(n+m+1-\alpha)} 2^{k_0 q \sum_{i=1}^{n} \frac{1}{p_i'}} 2^{kq \sum_{i=1}^{n} \frac{1}{p_i}} \\ &\leqslant C \sum_{k=k_0+4}^{\infty} 2^{(k_0 - k)(m+1-\alpha + \sum_{i=1}^{n} \frac{1}{p_i'})} \leqslant C.\end{aligned}$$

If $1 < p < \infty$, write

$$\begin{aligned}\| \mathscr{M}_N f \|_{\dot{K}_{\vec{p}}^{\alpha,q}(\mathbb{R}^n)}^q &\sum_{k=-\infty}^{\infty} 2^{k\alpha q} \left( \sum_{l=-\infty}^{\infty} |\lambda_l| \| (\mathscr{M}_N a_l) \chi_k \|_{L^{\vec{p}}(\mathbb{R}^n)} \right)^q \\ &\leqslant C \sum_{k=-\infty}^{\infty} 2^{k\alpha q} \left( \sum_{l=k-1}^{\infty} |\lambda_l| \| a_l \|_{L^{\vec{p}}(\mathbb{R}^n)} \right)^q \\ &\quad + C \sum_{k=-\infty}^{\infty} 2^{k\alpha q} \left( \sum_{l=-\infty}^{k-2} |\lambda_l| \| (\mathscr{M}_N a_l) \chi_k \|_{L^{\vec{p}}(\mathbb{R}^n)} \right)^q \\ &= III + IV.\end{aligned}$$

Using Hölder's inequality on mixed-norm spaces, we get

$$III \leqslant C \sum_{k=-\infty}^{\infty} 2^{k\alpha q} \left( \sum_{l=k-1}^{\infty} |\lambda_l| |B_l|^{\frac{\alpha}{n}} \right)^q$$

$$\leqslant C \sum_{k=-\infty}^{\infty} 2^{k\alpha q} \left( \sum_{l=k-1}^{\infty} |\lambda_l|^q |B_l|^{\frac{-\alpha q}{2n}} \right) \left( \sum_{l=k-1}^{\infty} |B_l|^{\frac{-\alpha q'}{2n}} \right)^{\frac{q}{q'}}$$

$$\leqslant C \sum_{k=-\infty}^{\infty} 2^{k\alpha q/2} \sum_{l=k-1}^{\infty} |\lambda_l|^p |B_l|^{\frac{-\alpha q}{2n}}$$

$$\leqslant C \sum_{l=-\infty}^{\infty} |\lambda_l|^p.$$

Now suppose $\alpha - \sum_{i=1}^{n} 1/p_i' < m+1$. As in the argument for $II$, we can obtain that

$$IV \leqslant C \sum_{k=-\infty}^{\infty} \left( \sum_{l=-\infty}^{k-2} |\lambda_l| 2^{l(m+1)-k(n+m+1)} \left( \frac{|B_k|}{|B_l|} \right)^{\frac{\alpha}{n}} \|\chi_{B_l}\|_{L^{\vec{p}'}(\mathbb{R}^n)} \|\chi_{B_k}\|_{L^{\vec{p}}(\mathbb{R}^n)} \right)^q$$

$$\leqslant C \sum_{k=-\infty}^{\infty} \left( \sum_{l=-\infty}^{k-2} |\lambda_l| 2^{(l-k)\left(m+1-\alpha+\sum_{i=1}^{n}\frac{1}{p_i'}\right)} \right)^q$$

$$\leqslant C \sum_{k=-\infty}^{\infty} \left( \sum_{l=-\infty}^{k-2} |\lambda_l|^p 2^{(l-k)\left(m+1-\alpha+\sum_{i=1}^{n}\frac{1}{p_i'}\right)\frac{q}{2}} \right) \left( \sum_{l=-\infty}^{k-2} 2^{(l-k)\left(m+1-\alpha+\sum_{i=1}^{n}\frac{1}{p_i'}\right)\frac{q'}{2}} \right)^{\frac{q}{q'}}$$

$$\leqslant C \sum_{l=-\infty}^{\infty} |\lambda_l|^p.$$

The proof is finished. $\square$

## 4  Molecular characterization of mixed Herz-Hardy spaces

In this section, we will obtain the molecular decomposition of mixed Herz-type Hardy spaces. We first give the notation of central $(\alpha, \vec{p}; s, \varepsilon)_l$-molecule.

**Definition 4.1.** Let $n - \sum_{i=1}^{n} 1/p_i \leqslant \alpha < \infty$, $0 < q < \infty$, $1 < \vec{p} < \infty$, and $s \geqslant [\alpha - n + \sum_{i=1}^{n} 1/p_i]$ be a non-negative integer. Set $\varepsilon > \max\{s/n, \alpha/n + 1/n\sum_{i=1}^{n} 1/p_i - 1\}$, $a = 1 - 1/n \sum_{i=1}^{n} 1/p_i - \alpha/n + \varepsilon$ and $b = 1 - 1/n \sum_{i=1}^{n} 1/p_i + \varepsilon$. A function $M_l \in L^{\vec{p}}(\mathbb{R}^n)$ with $l \in \mathbb{Z}$ (or $l \in \mathbb{N}$) is said to be a dyadic central $(\alpha, \vec{p}; s, \varepsilon)_l$-molecule (or dyadic central $(\alpha, \vec{p}; s, \varepsilon)_l$)-molecule of restricted type) if it satisfies

(d) $\|M_l\|_{L^{\vec{p}}(\mathbb{R}^n)} \leqslant 2^{-l\alpha}$.

(e) $\mathcal{R}_{\vec{p}}(M_l) = \|M_l\|_{L^{\vec{p}}(\mathbb{R}^n)}^{a/b} \||\cdot|^{nb} M_l(\cdot)\|_{L^{\vec{p}}(\mathbb{R}^n)}^{1-a/b} < \infty$.

(f) $\int_{\mathbb{R}^n} M_l(x) x^\beta dx = 0$, for any $\beta$ with $|\beta| \leqslant s$.

**Definition 4.2.** Let $n - \sum_{i=1}^{n} 1/p_i \leqslant \alpha < \infty$, $0 < q < \infty$, $1 < \vec{p} < \infty$, and $s \geqslant [\alpha - n + \sum_{i=1}^{n} 1/p_i]$ be a non-negative integer. Set $\varepsilon > \max\{s/n, \alpha/n + v - 1\}$, $a = 1 - 1/n \sum_{i=1}^{n} 1/p_i - \alpha/n + \varepsilon$ and $b = 1 - 1/n \sum_{i=1}^{n} 1/p_i + \varepsilon$.

(i) A function $M \in L^{\vec{p}}(\mathbb{R}^n)$ is said to be a central $(\alpha, \vec{p}; s, \varepsilon)$-molecule if it satisfies

(g) $\mathcal{R}_{\vec{p}}(M) = \|M\|_{L^{\vec{p}}(\mathbb{R}^n)}^{a/b} \||\cdot|^{nb} M(\cdot)\|_{L^{\vec{p}}(\mathbb{R}^n)}^{1-a/b} < \infty$.

(h) $\int_{\mathbb{R}^n} M(x) x^\beta dx = 0$, for any $\beta$ with $|\beta| \leqslant s$.

(ii) A function $M \in L^{\vec{p}}(\mathbb{R}^n)$ is said to be a central $(\alpha, \vec{p}; s, \varepsilon)$-molecule of restricted type if it satisfies (g),(h) in (i) and (g)' $\|M\|_{L^{\vec{p}}(\mathbb{R}^n)} \leqslant 1$.

The following lemma implies that the molecule is a generalization of atom.

**Lemma 4.1.** *Let $n - \sum_{i=1}^{n} 1/p_i \leqslant \alpha < \infty$, $0 < q < \infty$, $1 < \vec{p} < \infty$, and $s \geqslant [\alpha - n + \sum_{i=1}^{n} 1/p_i]$ be a non-negative integer. Set $\varepsilon > \max\{s/n, \alpha/n + 1/n \sum_{i=1}^{n} 1/p_i - 1\}$, $a = 1 - 1/n \sum_{i=1}^{n} 1/p_i - \alpha/n + \varepsilon$ and $b = 1 - 1/n \sum_{i=1}^{n} 1/p_i + \varepsilon$. If $M$ is a central $(\alpha, \vec{p})$-atom (or $(\alpha, \vec{p})$-atom of restricted type), then $M$ is also a central $(\alpha, \vec{p}; s, \varepsilon)$-molecule (or $(\alpha, \vec{p}; s, \varepsilon)$-molecule of restricted type) such that $\mathcal{R}_{\vec{p}}(M) \leqslant C$ with $C$ independent of $M$.*

*Proof.* We only need consider the case that $a$ is a $(\alpha, \vec{p})$-atom with support on a ball $B(0, r)$. A straightforward computation leads to that

$$\|M\|_{L^{\vec{p}}(\mathbb{R}^n)} \leqslant 2^{-\alpha l}$$

and

$$\|M\|_{L^{\vec{p}}(\mathbb{R}^n)}^{a/b} \left\| |\cdot|^{nb} M(\cdot) \right\|_{L^{\vec{p}}(\mathbb{R}^n)}^{1-a/b} \leqslant r^{nb(1-a/b)} \|M\|_{L^{\vec{p}}(\mathbb{R}^n)} \leqslant C r^\alpha r^{-\alpha} \leqslant C.$$

□

Now we give the molecular decomposition of mixed Herz-type Hardy spaces.

**Theorem 4.1.** *Let $n - \sum_{i=1}^{n} 1/p_i \leqslant \alpha < \infty$, $0 < q < \infty$, $1 < \vec{p} < \infty$, and $s \geqslant [\alpha - n + \sum_{i=1}^{n} 1/p_i]$ be a non-negative integer. Set $\varepsilon > \max\{s/n, \alpha/n + 1/n \sum_{i=1}^{n} 1/p_i - 1\}$, $a = 1 - 1/n \sum_{i=1}^{n} 1/p_i - \alpha/n + \varepsilon$ and $b = 1 - 1/n \sum_{i=1}^{n} 1/p_i + \varepsilon$. Then we have*
   *(i) $f \in H\dot{K}_{\vec{p}}^{\alpha,q}(\mathbb{R}^n)$ if and only if $f$ can be represented as*

$$f = \sum_{k=-\infty}^{\infty} \lambda_k M_k, \quad \text{in the sense of} \quad \mathcal{S}'(\mathbb{R}^n),$$

*where each $M_k$ is a dyadic central $(\alpha, \vec{p}; s, \varepsilon)_k$-molecule, and $\sum_{k=-\infty}^{\infty} |\lambda_k|^q < \infty$. Moreover,*

$$\|f\|_{H\dot{K}_{\vec{p}}^{\alpha,q}(\mathbb{R}^n)} \approx \inf \left( \sum_{k=-\infty}^{\infty} |\lambda_k|^q \right)^{1/q},$$

*where the infimum is taken over all above decompositions of $f$.*
   *(ii) $f \in HK_{\vec{p}}^{\alpha,q}(\mathbb{R}^n)$ if and only if*

$$f = \sum_{k=0}^{\infty} \lambda_k M_k, \quad \text{in the sense of} \quad \mathcal{S}'(\mathbb{R}^n),$$

*where each $M_k$ is a dyadic central $(\alpha, \vec{p}; s, \varepsilon)_k$-molecule of restricted type, and $\sum_{k=0}^{\infty} |\lambda_k|^q < \infty$. Moreover,*

$$\|f\|_{HK_{\vec{p}}^{\alpha,q}(\mathbb{R}^n)} \approx \inf \left( \sum_{k=0}^{\infty} |\lambda_k|^q \right)^{1/q},$$

*where the infimum is taken over all above decompositions of $f$.*

**Theorem 4.2.** *Let $n - \sum_{i=1}^{n} 1/p_i \leqslant \alpha < \infty$, $0 < q < \infty$, $1 < \vec{p} < \infty$, and $s \geqslant [\alpha - n + \sum_{i=1}^{n} 1/p_i]$ be a non-negative integer. Set $\varepsilon > \max\{s/n, \alpha/n + \frac{1}{n}\sum_{i=1}^{n} 1/p_i - 1\}$, $a = 1 - \frac{1}{n}\sum_{i=1}^{n} 1/p_i - \alpha/n + \varepsilon$ and $b = 1 - \frac{1}{n}\sum_{i=1}^{n} 1/p_i + \varepsilon$. Then we have*

(i) $f \in H\dot{K}_{\vec{p}}^{\alpha,q}(\mathbb{R}^n)$ if and only if $f$ can be represented as

$$f = \sum_{k=1}^{\infty} \lambda_k M_k, \quad \text{in the sense of} \quad \mathcal{S}'(\mathbb{R}^n),$$

where each $M_k$ is a central $(\alpha, \vec{p}; s, \varepsilon)$-molecule, and $\sum_{k=1}^{\infty} |\lambda_k|^q < \infty$. Moreover,

$$\|f\|_{H\dot{K}_{\vec{p}}^{\alpha,q}(\mathbb{R}^n)} \approx \inf \left( \sum_{k=1}^{\infty} |\lambda_k|^q \right)^{1/q},$$

where the infimum is taken over all above decompositions of $f$.

(ii) $f \in HK_{\vec{p}}^{\alpha,q}(\mathbb{R}^n)$ if and only if

$$f = \sum_{k=1}^{\infty} \lambda_k M_k, \quad \text{in the sense of} \quad \mathcal{S}'(\mathbb{R}^n),$$

where each $M_k$ is a central $(\alpha, \vec{p}; s, \varepsilon)$-molecule of restricted type, and $\sum_{k=1}^{\infty} |\lambda_k|^q < \infty$. Moreover,

$$\|f\|_{HK_{\vec{p}}^{\alpha,q}(\mathbb{R}^n)} \approx \inf \left( \sum_{k=1}^{\infty} |\lambda_k|^q \right)^{1/q},$$

where the infimum is taken over all above decompositions of $f$.

**Lemma 4.2.** *Let $n - \sum_{i=1}^{n} 1/p_i \leq \alpha < \infty$, $0 < q < \infty$, $1 < \vec{p} < \infty$, and $s \geq [\alpha - n + \sum_{i=1}^{n} 1/p_i]$ be a non-negative integer. Set $\varepsilon > \max\{s/n, \alpha/n + \frac{1}{n}\sum_{i=1}^{n} 1/p_i - 1\}$, $a = 1 - \frac{1}{n}\sum_{i=1}^{n} 1/p_i - \alpha/n + \varepsilon$ and $b = 1 - \frac{1}{n}\sum_{i=1}^{n} 1/p_i + \varepsilon$.*

*(i) If $0 < q \leq 1$, there exists a constant $C$ such that for any central $(\alpha, \vec{p}; s, \varepsilon)$ molecule (or $(\alpha, \vec{p}; s, \varepsilon)$-molecule of restricted type) $M$,*

$$\|M\|_{H\dot{K}_{\vec{p}}^{\alpha,q}(\mathbb{R}^n)} \leq C \left( \text{ or } \|M\|_{HK_{\vec{p}}^{\alpha,q}(\mathbb{R}^n)} \leq C \right).$$

*(ii) There exists a constant $C$ such that for any $l \in \mathbb{Z}$ (or $l \in \mathbb{N}$) and dyadic central $(\alpha, \vec{p}; s, \varepsilon)_l$-molecule (or dyadic central $(\alpha, \vec{p}; s, \varepsilon)_l$-molecule of restricted type) $M_l$,*

$$\|M_l\|_{H\dot{K}_{\vec{p}}^{\alpha,q}(\mathbb{R}^n)} \leq C \left( \text{ or } \|M_l\|_{HK_{\vec{p}}^{\alpha,q}(\mathbb{R}^n)} \leq C \right).$$

*Proof.* We just prove (i) for the homogeneous case, since the non-homogeneous case and the proof of (ii) can be proved similarly.

Let $M$ be a central $(\alpha, \vec{p}; s, \varepsilon)$-molecule. Write

$$\sigma = \|M\|_{L^{\vec{p}}(\mathbb{R}^n)}^{-1/\alpha}, \quad E_0 = \{x : |x| \leq \sigma\}$$

and

$$E_{k,\sigma} = \left\{ x : 2^{k-1}\sigma < |x| \leq 2^k \sigma \right\}, k \in \mathbb{Z}_+.$$

Set $B_{k,\sigma} = \{x : |x| \leq 2^k \sigma\}$ and denote by $\chi_{k,\sigma}$ the characteristic function of $E_{k,\sigma}$. Immediately,

we have
$$M(x) = \sum_{k=0}^{\infty} M(x)\chi_{k,\sigma}(x).$$

Let $M_k(x) = M(x)\chi_{k,\sigma}(x)$, $\mathcal{P}_m$ be the class of all real polynomials of degree $m$, and $P_{E_{k,\sigma}} M_k \in \mathcal{P}_m$ be the unique polynomial satisfying
$$\int_{E_{k,\sigma}} \left(M_k(x) - P_{E_{k,\sigma}} M_k(x)\right) x^{\beta} dx = 0, \quad |\beta| \leqslant s.$$

Set $Q_k(x) = \left(P_{E_{k,\sigma}} M_k\right)(x)\chi_{k,\sigma}(x)$. If we can prove that
  (a) there is a constant $C > 0$ and a sequences of numbers $\{\lambda_k\}_{k=0}^{\infty}$ such that
$$\sum_{k=0}^{\infty} |\lambda_k|^q < \infty, M_k - Q_k = \lambda_k a_k,$$
where each $a_k$ is a $(\alpha, \vec{p})$-atom;
  (b) $\sum_{k=0}^{\infty} Q_k$ has a $(\alpha, \vec{p})$-atom decomposition; then our desired conclusion can be deduced directly.

We first show (a). Without loss of generality, we suppose that $\mathcal{R}_{\vec{p}}(M) = 1$, which implies
$$\left\| |\cdot|^{nb} M(\cdot) \right\|_{L^{\vec{p}}(\mathbb{R}^n)} = \|M\|_{L^{\vec{p}}(\mathbb{R}^n)}^{-a/(b-a)} = \sigma^{na}.$$

Let $\{\varphi_l^k : |l| \leqslant s\} \subset \mathcal{P}_s(\mathbb{R}^n)$ such that
$$\left\langle \varphi_\mu^k, \varphi_\nu^k \right\rangle_{E_{k,\sigma}} = \frac{1}{|E_{k,\sigma}|} \int_{E_{k,\sigma}} \varphi_\mu^k(x)\varphi_\nu^k(x) dx = \delta_{\mu\nu}.$$

It is easy to deduce that
$$Q_k(x) = \sum_{|l| \leqslant s} \left\langle M_k, \varphi_l^k \right\rangle_{E_{k,\sigma}} \varphi_l^k(x), \text{ if } x \in E_{k,\sigma}$$

and
$$|Q_k(x)| \leqslant \left| \sum_{|l| \leqslant s} \frac{1}{|E_{k,\sigma}|} \int_{E_{k,\sigma}} M_k(x)\varphi_l^k(x) dx \varphi_l^k(x) \right| \leqslant \frac{C}{|E_{k,\sigma}|} \int_{E_{k,\sigma}} |M_k(x)| \, dx.$$

Thus for any $k \in \mathbb{Z}_+$, we have
$$\|M_k - Q_k\|_{L^{\vec{p}}(\mathbb{R}^n)}$$
$$\leqslant \|M_k\|_{L^{\vec{p}}(\mathbb{R}^n)} + \|Q_k\|_{L^{\vec{p}}(\mathbb{R}^n)}$$
$$\leqslant \|M_k\|_{L^{\vec{p}}(\mathbb{R}^n)} + \frac{C}{|E_{k,\sigma}|} \|M_k\|_{L^{\vec{p}}(\mathbb{R}^n)} \|\chi_{E_{k,\sigma}}\|_{L^{\vec{p}'}(\mathbb{R}^n)} \|\chi_{E_{k,\sigma}}\|_{L^{\vec{p}}(\mathbb{R}^n)}$$
$$\leqslant \|M_k\|_{L^{\vec{p}}(\mathbb{R}^n)} + C \|M_k\|_{L^{\vec{p}}(\mathbb{R}^n)}$$
$$\leqslant C \|M_k\|_{L^{\vec{p}}(\mathbb{R}^n)}$$

$$\leqslant C \left\| |\cdot|^{nb} M(\cdot) \right\|_{L^{\vec{p}}(\mathbb{R}^n)} \left| 2^k \sigma \right|^{-nb}$$

$$= C \left| 2^k \sigma \right|^{-nb} \sigma^{na} = C 2^{-kna} |B_{k,\sigma}|^{-\alpha}.$$

We see that $M_k - Q_k = \lambda_k a_k$, with $\lambda_k = C 2^{-kn\alpha}$ and $a_k$ a central $(\alpha, \vec{p})$-atom supported in $B_{k,\sigma}$. Moreover, one can easily get $\sum_{k=0}^{\infty} |\lambda_k|^q < \infty$.

Next we will show (b). Let $\{\psi_l^k : |l| \leqslant s\} \subset \mathcal{P}_s(\mathbb{R}^n)$ be the dual basis of $\{x^\alpha : |\alpha| \leqslant s\}$ with respect to the weight $1/|E_{k,\sigma}|$ on $E_{k,\sigma}$, that is

$$\left\langle \psi_l^k, x^\alpha \right\rangle = \frac{1}{|E_{k,\sigma}|} \int_{E_{k,\sigma}} \psi_l^k(x) x^\alpha dx = \delta_{l\alpha}.$$

Setting $\varphi_l^k(x) = \sum_{|\nu| \leqslant s} \beta_{l\nu}^k x^\nu$ and $\psi_l^k(x) = \sum_{|\nu| \leqslant s} \tau_{\nu l}^k \varphi_\nu^k(x)$, we have

$$\tau_{\nu l}^k = \left\langle \psi_l^k, \varphi_\nu^k \right\rangle = \sum_{|\gamma| \leqslant s} \beta_{\nu\gamma}^k \left\langle \psi_l^k, x^\gamma \right\rangle = \sum_{|\gamma| \leqslant s} \beta_{\nu\gamma}^k \delta_{l\gamma} = \beta_{\nu l}^k.$$

Thus

$$\psi_l^k(x) = \sum_{|\nu| \leqslant s} \beta_{\nu l}^k \varphi_\nu^k(x). \tag{6}$$

By using (6), we deduce that for $x \in E_{k,\sigma}$,

$$\left\langle M_k, \varphi_l^k \right\rangle_{E_{k,\sigma}} \varphi_l^k(x) = \left\langle M_k, \sum_{|\nu| \leqslant s} \beta_{\nu l}^k x^\nu \right\rangle_{E_{k,\sigma}} \varphi_l^k(x) = \sum_{|\nu| \leqslant s} \left\langle M_k, x^\nu \right\rangle_{E_{k,\sigma}} \beta_{\nu l}^k \varphi_l^k(x),$$

which together with $Q_k(x) = \sum_{|l| \leqslant s} \left\langle M_k, \varphi_l^k \right\rangle_{E_{k,\sigma}} \varphi_l^k(x)$, implies that

$$Q_k(x) = \sum_{|l| \leqslant s} \left\langle M_k, x^l \right\rangle_{E_{k,\sigma}} \psi_l^k(x), \text{ if } x \in E_{k,\sigma}.$$

Next we will show that there is a constant $C > 0$ such that

$$\left| \psi_l^k(x) \right| \leqslant C \left( 2^{k-1} \sigma \right)^{-|l|}. \tag{7}$$

Set $E = \{x \in \mathbb{R}^n : 1 \leqslant |x| \leqslant 2\}$, $F = \{x \in \mathbb{R}^n : |x| \leqslant 1\}$, $\{e_l : |l| \leqslant s\} \subset \mathcal{P}_s(\mathbb{R}^n)$ satisfying $\frac{1}{|E|} \int_E e_l(x) x^\alpha dx = \delta_{l\alpha}$, and $\{\tilde{e}_l : |l| \leqslant s\} \subset \mathcal{P}_s(\mathbb{R}^n)$ satisfying $\frac{1}{|F|} \int_F \tilde{e}_l(x) x^\alpha dx = \delta_{l\alpha}$.

Noting that

$$\delta_{l\alpha} = \frac{1}{|E_{k,\sigma}|} \int_{E_{k,\sigma}} \psi_l^k(x) x^\alpha dx = \frac{1}{|E|} \int_E \left( 2^{k-1} \sigma \right)^{|\alpha|} \psi_l^k \left( 2^{k-1} \sigma y \right) y^\alpha dy,$$

we get $e_l(y) = \left( 2^{k-1} \sigma \right)^{|l|} \psi_l^k \left( 2^{k-1} \sigma y \right)$. This in turn leads to that

$$\psi_l^k(y) = \left( 2^{k-1} \sigma \right)^{-|l|} e_l \left( \frac{x}{2^{k-1} \sigma} \right), x \in E_{k,\sigma}.$$

Similarly,
$$\psi_l^k(x) = \left(2^{k-1}\sigma\right)^{-|l|} \tilde{e}_l\left(\frac{x}{\sigma}\right), \quad x \in F.$$

Taking $C = \sup_{l:|l|\leqslant s}\left\{\|e_l\|_{L^\infty(\mathbb{R}^n)}, \|\tilde{e}_l\|_{L^\infty(\mathbb{R}^n)}\right\}$. This inequality (7) then follows directly.

We can now complete the proof of (b).

Setting
$$N_l^k = \sum_{j=k}^{\infty} |E_{j,\sigma}| \left\langle M_j, x^l \right\rangle_{E_{j,\sigma}}, \quad k \in \mathbb{N},$$

it is readily to see that
$$N_l^0 = \sum_{j=0}^{\infty} |E_{j,\sigma}| \left\langle M_j, x^l \right\rangle_{E_{j,\sigma}} = \sum_{j=0}^{\infty} \int_{E_{j,\sigma}} M(x) x^l dx = \int_{\mathbb{R}^n} M(x) x^l dx = 0,$$

and for $k \in \mathbb{Z}_+$, there exists $E_\sigma \subset E_{j,\sigma}$ such that $|E_\sigma| = \min\{1, |E_{j,\sigma}|\}$. Therefore,

$$\left|N_l^k\right| \leqslant \sum_{j=k}^{\infty} \int_{E_{j,\sigma}} \left|M_j(x) x^l\right| dx$$

$$\leqslant C \sum_{j=k}^{\infty} \left\||M_j(\cdot)| \cdot |\cdot|^l\right\|_{L^{\vec{p}}(\mathbb{R}^n)} \|\chi_{E_{j,\sigma}}\|_{L^{\vec{p}'}(\mathbb{R}^n)}$$

$$\leqslant C \sum_{j=k}^{\infty} (2^j \sigma)^{|l|-nb} \left\||\cdot|^{nb} M_j(\cdot)\right\|_{L^{\vec{p}}(\mathbb{R}^n)} \left(\frac{|E_{j,\sigma}|}{|E_\sigma|}\right)^{\frac{1}{n}\sum_{i=1}^{n}\frac{1}{p_i'}} \|\chi_{E_\sigma}\|_{L^{\vec{p}'}(\mathbb{R}^n)}$$

$$\leqslant C \sum_{j=k}^{\infty} \sigma^{|l|+na-nb+\sum_{i=1}^{n}\frac{1}{p_i'}} 2^{j\left(|l|-nb+\sum_{i=1}^{n}\frac{1}{p_i'}\right)}$$

$$\leqslant C\sigma^{|l|-\alpha+\sum_{i=1}^{n}\frac{1}{p_i'}} 2^{k\left(|l|-nb+\sum_{i=1}^{n}\frac{1}{p_i'}\right)}.$$

This, together with (7) shows that
$$|E_{k,\sigma}|^{-1} \left|N_l^k \psi_l^k(x) \chi_{k,\sigma}(x)\right| \leqslant C\sigma^{\sum_{i=1}^{n}\frac{1}{p_i'}-n-\alpha} 2^{-kn\left(b+1-\frac{1}{n}\sum_{i=1}^{n}\frac{1}{p_i'}\right)} \to 0, \quad \text{if} \quad k \to \infty. \qquad (8)$$

By using (8) and Abel transform, it yields
$$\sum_{k=0}^{\infty} Q_k(x) = \sum_{|l|\leqslant s} \sum_{k=0}^{\infty} \left(\sum_{j=0}^{k} |E_{j,\sigma}| \left\langle M_j, x^l \right\rangle_{E_{j,\sigma}}\right)$$
$$\times \left\{|E_{k,\sigma}|^{-1} \psi_l^k(x) \chi_{k,\sigma}(x) - |E_{k+1,\sigma}|^{-1} \psi_l^{k+1}(x) \chi_{k+1,\sigma}(x)\right\}$$
$$= \sum_{|l|\leqslant s} \sum_{k=0}^{\infty} \left(-N_l^{k+1}\right) \left\{|E_{k,\sigma}|^{-1} \psi_l^k(x) \chi_{k,\sigma}(x) - |E_{k+1,\sigma}|^{-1} \psi_l^{k+1}(x) \chi_{k+1,\sigma}(x)\right\}.$$

Meanwhile,
$$\left|N_l^{k+1} \left\{|E_{k,\sigma}|^{-1} \psi_l^k(x) \chi_{k,\sigma}(x) - |E_{k+1,\sigma}|^{-1} \psi_l^{k+1}(x) \chi_{k+1,\sigma}(x)\right\}\right|$$

$$\leqslant C \left|N_l^{k+1}\right| |E_{k+1,\sigma}|^{-1} \left|\psi_l^{k+1}(x)\right|$$

$$\leqslant C 2^{-kna} |E_{k+1,\sigma}|^{\frac{1}{n}\sum_{i=1}^{n}\frac{1}{p_i'}-1-\alpha/n}.$$

Set $\lambda_{lk} = C 2^{-kna}$ and

$$a_{lk} = \lambda_{lk}^{-1}\left(-N_l^{k+1}\right)\left\{|E_{k,\sigma}|^{-1}\psi_l^k(x)\chi_{k,\sigma}(x) - |E_{k+1,\sigma}|^{-1}\psi_l^{k+1}(x)\chi_{k+1,\sigma}(x)\right\}.$$

Then

$$\sum_{k=0}^{\infty} Q_k(x) = \sum_{|l|\leqslant s}\sum_{k=0}^{\infty} \lambda_{lk} a_{lk}$$

with $a_{lk}$ a $(\alpha, \vec{p})$-atom, and $\sum_{|l|\leqslant s}\sum_{k=0}^{\infty} |\lambda_{lk}|^q < \infty$. The conclusion (b) then holds. $\square$

## 5 Application

In this section, we will give two applications for atom decomposition and molecular decomposition of mixed Herz-Hardy spaces. The first application is to build a boundedness criterion for certain sublinear operators on mixed Herz-Hardy spaces, and the other application is to give the dual spaces of mixed Herz-Hardy spaces.

### 5.1 Boundedness for a class of sublinear operators on mixed Herz-Hardy spaces

This subsection is devoted to building a boundedness criterion for certain sublinear operators from $H\dot{K}_{\vec{p}}^{\alpha,q}(\mathbb{R}^n)$ to $\dot{K}_{\vec{p}}^{\alpha,q}(\mathbb{R}^n)$ ( or from $HK_{\vec{p}}^{\alpha,q}(\mathbb{R}^n)$ to $K_{\vec{p}}^{\alpha,q}(\mathbb{R}^n)$). The result can be stated as follows.

**Theorem 5.1.** *Let $\sum_{i=1}^{n} 1/p_i' \leqslant \alpha < \infty, 0 < q < \infty, 1 < \vec{p} < \infty$ and the integer $s = [\alpha - \sum_{i=1}^{n} 1/p_i']$. Suppose $T$ is a sublinear operator satisfying:*

*(i) $T$ is bounded on $L^{\vec{p}}(\mathbb{R}^n)$;*

*(ii) there exists a constant $\delta > 0$ such that $s + \delta > \alpha - \sum_{i=1}^{n} 1/p_i'$, and for any compact support function $f$ with*

$$\int_{\mathbb{R}^n} f(x) x^\beta dx = 0, \quad |\beta| \leqslant s,$$

*$Tf$ satisfies the size condition*

$$|Tf(x)| \leqslant C(\operatorname{diam}(\operatorname{supp} f))^{s+\delta}|x|^{-(n+s+\delta)}\|f\|_{L^1(\mathbb{R}^n)}, \quad if \quad \operatorname{dist}(x, \operatorname{supp} f) \geqslant |x|/2.$$

*Then $T$ can be extended to be a bounded operator from $H\dot{K}_{\vec{p}}^{\alpha,q}(\mathbb{R}^n)$ to $\dot{K}_{\vec{p}}^{\alpha,q}(\mathbb{R}^n)$ ( or bounded from $HK_{\vec{p}}^{\alpha,q}(\mathbb{R}^n)$ to $K_{\vec{p}}^{\alpha,q}(\mathbb{R}^n)$).*

*Proof.* It suffices to prove the homogeneous case. Suppose $f \in H\dot{K}_{\vec{p}}^{\alpha,q}(\mathbb{R}^n)$. By Theorem 3.1, we may rewrite $f$ as $f = \sum_{j=-\infty}^{\infty} \lambda_j b_j$ in the sense of $\mathcal{S}'(\mathbb{R}^n)$, where each $b_j$ is a central $(\alpha, \vec{p})$-atom with support contained in $B_j$ and

$$\|f\|_{H\dot{K}_{\vec{p}}^{\alpha,q}(\mathbb{R}^n)} \approx \inf \left( \sum_{j=-\infty}^{\infty} |\lambda_j|^q \right)^{1/q}.$$

Then, we get

$$\|Tf\|_{\dot{K}_{\vec{p}}^{\alpha,q}(\mathbb{R}^n)}^q = \sum_{k=-\infty}^{\infty} 2^{k\alpha q} \|(Tf)\chi_k\|_{L^{\vec{p}}(\mathbb{R}^n)}^q$$

$$\leqslant C \left( \sum_{k=-\infty}^{\infty} 2^{k\alpha q} \left( \sum_{j=-\infty}^{k-2} |\lambda_j| \|(Tb_j)\chi_k\|_{L^{\vec{p}}(\mathbb{R}^n)} \right)^q \right.$$

$$\left. + \sum_{k=-\infty}^{\infty} 2^{k\alpha q} \left( \sum_{j=k-1}^{\infty} |\lambda_j| \|(Tb_j)\chi_k\|_{L^{\vec{p}}(\mathbb{R}^n)} \right)^q \right)$$

$$= C(I_1 + I_2).$$

Let us first estimate $I_1$. By the size condition of $T$ and Hölder's inequality on mixed-norm spaces, it yields

$$|Tb_j(x)| \leqslant C \left( \operatorname{diam}(\operatorname{supp} b_j) \right)^{s+\delta} |x|^{-(n+s+\delta)} \|b_j\|_{L^1(\mathbb{R}^n)}$$

$$\leqslant C|x|^{-(n+s+\delta)} 2^{j(s+\delta)} \int_{B_j} |b_j(y)| \, dy$$

$$\leqslant C|x|^{-(n+s+\delta)} 2^{j(s+\delta)} \|b_j\|_{L^{\vec{p}}(\mathbb{R}^n)} \|\chi_{B_j}\|_{L^{\vec{p}'}(\mathbb{R}^n)}$$

$$\leqslant C 2^{j(s+\delta-\alpha)-k(s+\delta+n)} \|\chi_{B_j}\|_{L^{\vec{p}'}(\mathbb{R}^n)}.$$

As a consequence,

$$\|(Tb_j)\chi_k\|_{L^{\vec{p}}(\mathbb{R}^n)} \leqslant C 2^{j(s+\delta-\alpha)-k(s+\delta+n)} \|\chi_{B_j}\|_{L^{\vec{p}'}(\mathbb{R}^n)} \|\chi_{B_k}\|_{L^{\vec{p}}(\mathbb{R}^n)}$$

$$= C 2^{j(s+\delta-\alpha)-k(s+\delta+n)} 2^{j \sum_{i=1}^{n} \frac{1}{p_i'}} 2^{k\left(n - \sum_{i=1}^{n} \frac{1}{p_i'}\right)}$$

$$\leqslant C 2^{(j-k)\left(s+\delta+\sum_{i=1}^{n} \frac{1}{p_i'}\right) - j\alpha}.$$

Therefore, when $0 < q \leqslant 1$, by $\sum_{i=1}^{n} \frac{1}{p_i'} \leqslant \alpha < s + \delta + \sum_{i=1}^{n} \frac{1}{p_i'}$, we get

$$I_1 = \sum_{k=-\infty}^{\infty} 2^{k\alpha q} \left( \sum_{j=-\infty}^{k-2} |\lambda_j| \|(Tb_j)\chi_k\|_{L^{\vec{p}}(\mathbb{R}^n)} \right)^q$$

$$\leqslant C \sum_{k=-\infty}^{\infty} 2^{k\alpha q} \left( \sum_{j=-\infty}^{k-2} |\lambda_j|^q 2^{\left((j-k)\left(s+\delta+\sum_{i=1}^{n} \frac{1}{p_i'}\right) - j\alpha\right)q} \right)$$

$$= C \sum_{j=-\infty}^{\infty} |\lambda_j|^q \sum_{k=j+2}^{\infty} 2^{(j-k)\left(s+\delta+\sum_{i=1}^{n}\frac{1}{p_i'}-\alpha\right)q} \leqslant C \sum_{j=-\infty}^{\infty} |\lambda_j|^q.$$

For $1 < q < \infty$, noting that $\sum_{i=1}^{n} 1/p_i' \leqslant \alpha < s + \delta + \sum_{i=1}^{n} 1/p_i'$, by the estimate of $|Tb_j(x)|$ and Hölder's inequality, we have

$$I_1 \leqslant C \sum_{k=-\infty}^{\infty} 2^{k\alpha q} \left( \sum_{j=-\infty}^{k-2} |\lambda_j| 2^{(j-k)\left(s+\delta+\sum_{i=1}^{n}\frac{1}{p_i'}\right)-j\alpha} \right)^q$$

$$\leqslant C \sum_{k=-\infty}^{\infty} \left( \sum_{j=-\infty}^{k-2} |\lambda_j|^q 2^{(j-k)\left(s+\delta+\sum_{i=1}^{n}\frac{1}{p_i'}-\alpha\right)\frac{q}{2}} \right)$$

$$\times \left( \sum_{j=-\infty}^{k-2} 2^{(j-k)\left(s+\delta+\sum_{i=1}^{n}\frac{1}{p_i'}-\alpha\right)\frac{q'}{2}} \right)^{\frac{q}{q'}}$$

$$\leqslant C \sum_{k=-\infty}^{\infty} \left( \sum_{j=-\infty}^{k-2} |\lambda_j|^q 2^{(j-k)\left(s+\delta+\sum_{i=1}^{n}\frac{1}{p_i'}-\alpha\right)\frac{q}{2}} \right)$$

$$\leqslant C \sum_{j=-\infty}^{\infty} |\lambda_j|^q \sum_{k=j+2}^{\infty} 2^{(j-k)\left(s+\delta+\sum_{i=1}^{n}\frac{1}{p_i'}-\alpha\right)\frac{q}{2}}$$

$$\leqslant C \sum_{j=-\infty}^{\infty} |\lambda_j|^q.$$

Let us now estimate $I_2$. When $0 < q \leqslant 1$, by the $L^{\vec{p}}(\mathbb{R}^n)$ boundedness of $T$, there holds

$$I_2 = \sum_{k=-\infty}^{\infty} 2^{k\alpha q} \left( \sum_{j=k-1}^{\infty} |\lambda_j| \, \|(Tb_j)\chi_k\|_{L^{\vec{p}}(\mathbb{R}^n)} \right)^q$$

$$\leqslant C \sum_{k=-\infty}^{\infty} 2^{k\alpha q} \left( \sum_{j=k-1}^{\infty} |\lambda_j|^q \|b_j\|_{L^{\vec{p}}(\mathbb{R}^n)}^q \right)$$

$$\leqslant C \sum_{k=-\infty}^{\infty} 2^{k\alpha q} \left( \sum_{j=k-1}^{\infty} |\lambda_j|^q |B_j|^{\frac{-\alpha q}{n}} \right)$$

$$= C \sum_{j=-\infty}^{\infty} |\lambda_j|^q \sum_{k=-\infty}^{j+1} 2^{(k-j)\alpha q} \leqslant C \sum_{j=-\infty}^{\infty} |\lambda_j|^q.$$

When $1 < q < \infty$, again by the $L^{\vec{p}}(\mathbb{R}^n)$ boundedness of $T$ and Hölder's inequality, we have

$$I_2 \leqslant C \sum_{k=-\infty}^{\infty} 2^{k\alpha q} \left( \sum_{j=k-1}^{\infty} |\lambda_j|^q \|(Tb_j)\chi_k\|_{L^{\vec{p}}(\mathbb{R}^n)}^{\frac{q}{2}} \right) \times \left( \sum_{j=k-1}^{\infty} \|(Tb_j)\chi_k\|_{L^{\vec{p}}(\mathbb{R}^n)}^{\frac{q'}{2}} \right)^{\frac{q}{q'}}$$

$$\leqslant C \sum_{k=-\infty}^{\infty} 2^{k\alpha q} \left( \sum_{j=k-1}^{\infty} |\lambda_j|^q \|b_j\|_{L^{\vec{p}}(\mathbb{R}^n)}^{\frac{q}{2}} \right) \left( \sum_{j=k-1}^{\infty} \|b_j\|_{L^{\vec{p}}(\mathbb{R}^n)}^{\frac{q'}{2}} \right)^{\frac{q}{q'}}$$

$$\leqslant C \sum_{k=-\infty}^{\infty} 2^{k\alpha q} \left( \sum_{j=k-1}^{\infty} |\lambda_j|^q |B_j|^{\frac{-\alpha q}{2n}} \right) \left( \sum_{j=k-1}^{\infty} |B_j|^{\frac{-\alpha q'}{2n}} \right)^{\frac{q}{q'}}$$

$$\leqslant C \sum_{k=-\infty}^{\infty} 2^{k\alpha q/2} \left( \sum_{j=k-1}^{\infty} |\lambda_j|^q |B_j|^{\frac{-\alpha q}{2n}} \right)$$

$$= C \sum_{j=-\infty}^{\infty} |\lambda_j|^q \sum_{k=-\infty}^{j+1} 2^{\frac{(k-j)\alpha q}{2}} \leqslant C \sum_{j=-\infty}^{\infty} |\lambda_j|^q.$$

Combining all the estimates above, we arrive at

$$\|Tf\|_{\dot{K}_{\vec{p}}^{\alpha,q}(\mathbb{R}^n)} \leqslant C\|f\|_{H\dot{K}_{\vec{p}}^{\alpha,q}(\mathbb{R}^n)}.$$

Thus, the proof is completed. $\square$

Now we give some concrete operators satisfying the conditions of Theorem 5.1.

The classical Calderón-Zygmund operator $T$ is a initially $L^2(\mathbb{R}^n)$ bounded operator with the associated standard kernel $K$, that is to say, functions $K(x,y)$ defined on $\mathbb{R}^n \times \mathbb{R}^n \setminus \{(x,x) : x \in \mathbb{R}^n\}$ satisfying

(i) the size condition:
$$|K(x,y)| \leqslant \frac{A}{|x-y|^n},$$
for some constant $A > 0$;

(ii) the regularity conditions: for some $\delta > 0$,

$$\left|K(x,y) - K\left(x',y\right)\right| \leqslant \frac{A|x-x'|^{\delta}}{(|x-y|+|x'-y|)^{n+\delta}}$$

holds whenever $|x-x'| \leqslant \frac{1}{2} \max(|x-y|, |x'-y|)$ and

$$\left|K(x,y) - K\left(x,y'\right)\right| \leqslant \frac{A|y-y'|^{\delta}}{(|x-y|+|x-y'|)^{n+\delta}}$$

holds whenever $|y-y'| \leqslant \frac{1}{2} \max(|x-y|, |x-y'|)$.

Then $T$ can be represented as

$$Tf(x) = \int_{\mathbb{R}^n} K(x,y)f(y)dy, \quad x \notin \operatorname{supp} f,$$

which is obvious a $L^{\vec{p}}(\mathbb{R}^n)$-bounded operator.

**Theorem 5.2.** *Suppose that $T$ is an above Calderón-Zygmund operator, and that $0 < \delta \leqslant 1$ is the constant associated with the standard kernel $K$, then for $\sum_{i=1}^{n} 1/p_i' \leqslant \alpha < \sum_{i=1}^{n} 1/p_i' + \delta$ and $0 < q < \infty$, $T$ is bounded from $H\dot{K}_{\vec{p}}^{\alpha,q}(\mathbb{R}^n)$ to $\dot{K}_{\vec{p}}^{\alpha,q}(\mathbb{R}^n)$.*

*Proof.* Noting that $\sum_{i=1}^{n} 1/p_i' \leq \alpha < \sum_{i=1}^{n} 1/p_i' + \delta$ implies that $s = [\alpha - \sum_{i=1}^{n} 1/p_i'] = 0$, the operator $T$ satisfies the conditions of Theorem 5.1 with $s = 0$. The desired conclusion follows directly. $\square$

**Theorem 5.3.** *For any central $(\alpha, \vec{p})$-atom $f$, let*

$$Tf(x) = \int_{\mathbb{R}^n} K(x,y) f(y) dy, \quad x \notin \operatorname{supp} f$$

*satisfying $\int_{\mathbb{R}^n} Tf(x) dx = 0$ be a bounded operator on $L^{\vec{p}}(\mathbb{R}^n)$ for some $1 < \vec{p} < \infty$, and the kernel $K$ satisfies that there are constants $C' > 0$ and $0 < \delta \leq 1$ such that*

$$|K(x,y) - K(x,0)| \leq C' \frac{|y|^\delta}{|x-y|^{n+\delta}}, \quad |x| \geq 2|y|.$$

*Then for any $\alpha$ and $q$ with $\sum_{i=1}^{n} 1/p_i' \leq \alpha < \sum_{i=1}^{n} 1/p_i' + \delta$ and $0 < q < \infty$, there exists a constant $C$ such that $\|Tf\|_{H\dot{K}^{\alpha,q}_{\vec{p}}(\mathbb{R}^n)} \leq C$ ( or $\|Tf\|_{HK^{\alpha,q}_{\vec{p}}(\mathbb{R}^n)} \leq C$).*

*Proof.* We only prove homogeneous case. Let $f$ be a central $(\alpha, \vec{p})$-atom supporting in $B(0,r)(r > 0)$. It suffices to show $Tf$ is a central $(\alpha, \vec{p}; 0, \varepsilon)$-molecule for some $1 + \delta/n - 1/n \sum_{i=1}^{n} 1/p_i' \geq \varepsilon > \alpha/n - 1/n \sum_{i=1}^{n} 1/p_i'$. Let $a = 1/n \sum_{i=1}^{n} 1/p_i' - \alpha/n + \varepsilon$, $b = 1/n \sum_{i=1}^{n} 1/p_i' + \varepsilon$. Next, we will verify the size condition for molecules, that is

$$\mathcal{R}_{\vec{p}}(Tf) = \|Tf\|_{L^{\vec{p}}(\mathbb{R}^n)}^{a/b} \left\| |\cdot|^{nb}(Tf)(\cdot) \right\|_{L^{\vec{p}}(\mathbb{R}^n)}^{1-a/b} \leq C,$$

with $C$ independent of $f$. To do this, we first estimate $\left\| |\cdot|^{nb}(Tf)(\cdot) \right\|_{L^{\vec{p}}(\mathbb{R}^n)}$ In fact, we have

$$\left\| |\cdot|^{nb}(Tf)(\cdot) \right\|_{L^{\vec{p}}(|\cdot| \leq 2r)} \leq Cr^{nb} \|Tf\|_{L^{\vec{p}}(\mathbb{R}^n)} \leq Cr^{nb-\alpha}.$$

Moreover, the vanishing moment of $f$ and the regularity of $K$ give us that for $x$ with $|x| > 2r$,

$$\begin{aligned}
|Tf(x)| &= \left| \int_{\mathbb{R}^n} K(x,y) f(y) dy \right| \\
&= \left| \int_{\mathbb{R}^n} (K(x,y) - K(x,0)) f(y) dy \right| \\
&\leq C \int_{\mathbb{R}^n} \frac{|y|^\delta}{|x-y|^{n+\delta}} |f(y)| dy \\
&\leq Cr^{n+\delta} |x|^{-(n+\delta)} \frac{1}{|B_{0,r}|} \int_{B_{0,r}} |f(y)| dy \\
&\leq Cr^{n+\delta} |x|^{-(n+\delta)} Mf(x).
\end{aligned}$$

Therefore, by the fact $nb - n - \delta \leq 0$, we further get

$$\begin{aligned}
\left\| |\cdot|^{nb}(Tf)(\cdot) \right\|_{L^{\vec{p}}(|\cdot|>2r)} &\leq Cr^{n+\delta} \left\| |\cdot|^{nb-n-\delta} Mf(\cdot) \right\|_{L^{\vec{p}}(|\cdot|>2r)} \\
&\leq Cr^{n+\delta+nb-n-\delta} \|Mf\|_{L^{\vec{p}}(\mathbb{R}^n)}
\end{aligned}$$

$$\leqslant Cr^{nb}\|f\|_{L^{\vec{p}}(\mathbb{R}^n)} \leqslant Cr^{nb-\alpha}.$$

Consequently,

$$\begin{aligned}\mathcal{R}_{\vec{p}}(Tf) &= \|Tf\|_{L^{\vec{p}}(\mathbb{R}^n)}^{a/b} \left\||\cdot|^{nb}(Tf)(\cdot)\right\|_{L^{\vec{p}}(\mathbb{R}^n)}^{1-a/b} \\ &\leqslant C\|f\|_{L^{\vec{p}}(\mathbb{R}^n)}^{a/b} r^{(nb-\alpha)(1-a/b)} \\ &\leqslant Cr^{-\alpha a/b + (nb-\alpha)(1-a/b)} \leqslant C.\end{aligned}$$

□

## 5.2 Dual spaces of mixed homogeneous Herz-Hardy spaces

This subsection will investigate the dual spaces of mixed homogeneous Herz-Hardy spaces $H\dot{K}_{\vec{p}}^{\alpha,q}(\mathbb{R}^n)$ by using the atom decomposition.

The duality theory plays a key role in functional analysis and harmonic analysis. As is well known, the dual space of the classical Hardy space $H^1(\mathbb{R}^n)$ is $BMO(\mathbb{R}^n)$. In 2019, Huang [26] investigated the dual space of mixed Hardy spaces and pointed out that they are Campanato-type spaces, which can reduce to $BMO(\mathbb{R}^n)$ in some special cases. Recently, Wei [24] proved that the dual spaces of some particular mixed atomic Hardy spaces are mixed central BMO spaces. Note that the mixed atomic Hardy spaces studied in [24] are special cases of the mixed homogeneous Herz-Hardy spaces by in view of Theorem 3.1, it is natural for us to consider the dual space of mixed homogeneous Herz spaces. In this subsection, we will prove that the dual spaces of mixed homogeneous Herz spaces are mixed central Campanato spaces. To this end, we first give the definition of the mixed central Campanato space $\mathcal{CL}_{\vec{p},\alpha,s}(\mathbb{R}^n)$.

For any $s \in \mathbb{Z}_+$, $\mathcal{P}_s(\mathbb{R}^n)$ denotes the set of all polynomials on $\mathbb{R}^n$ with degree not greater than $s$. For any ball $B$ and any locally integrable function $g$ on $\mathbb{R}^n$, we use $P_B^s g$ to denote the minimizing polynomial of $g$ with degree not greater than $s$, which means that $P_B^s g$ is the unique polynomial $f \in \mathcal{P}_s(\mathbb{R}^n)$ such that, for any $h \in \mathcal{P}_s(\mathbb{R}^n)$,

$$\int_B [g(x) - f(x)]h(x)\mathrm{d}x = 0.$$

**Definition 5.1.** *Let* $n(1-1/n\sum_{i=1}^m 1/p_i) \leqslant \alpha < \infty$, *non-negative integer* $s \geqslant [\alpha - n(1 - 1/n\sum_{i=1}^m 1/p_i)]$, *and* $1 < \vec{p} < \infty$. *Then the mixed central Campanato space* $\mathcal{CL}_{\vec{p},\alpha,s}(\mathbb{R}^n)$ *is defined by*

$$\|f\|_{\mathcal{CL}_{\vec{p},\alpha,s}(\mathbb{R}^n)} = \sup_{r>0} |B(0,r)|^{\frac{-\alpha}{n}} \left\|\left(f - P_{B(0,r)}^s f\right)\chi_{B(0,r)}\right\|_{L^{\vec{p}}(\mathbb{R}^n)} < \infty.$$

**Theorem 5.4.** *Let* $1 < \vec{p} < \infty$, *Then the dual space of the mixed Herz-Hardy space* $H\dot{K}_{\vec{p}}^{\alpha,q}(\mathbb{R}^n)$ *is the mixed central bounded mean oscillation space* $\mathcal{CL}_{\vec{p},\alpha,s}(\mathbb{R}^n)$, *that is*

$$\left(H\dot{K}_{\vec{p}}^{\alpha,q}(\mathbb{R}^n)\right)^* = \mathcal{CL}_{\vec{p'},\alpha,s}(\mathbb{R}^n).$$

*Proof.* We divide our proof into two steps.

(i) $\mathcal{CL}_{\vec{p'},\alpha,s}(\mathbb{R}^n) \subseteq \left(H\dot{K}_{\vec{p}}^{\alpha,q}(\mathbb{R}^n)\right)^*$, that is to say, any $f \in \mathcal{CL}^{\vec{p'},\alpha,s}(\mathbb{R}^n)$ can be seen as a bounded linear functional on $H\dot{K}_{\vec{p}}^{\alpha,q}(\mathbb{R}^n)$.

Without loss of generality, we just prove above conclusion for any central $(\alpha, \vec{p})$ atom $a_j$. Assume $a_j$ is supported in $B_j$, $j \in \mathbb{Z}$. Then

$$\left| \int_{B_j} g(x) a(x) dx \right| = \left| \int_{B_j} \left( g(x) - P^s_{B_j} g(x) \right) a(x) dx \right|$$

$$\leqslant \left\| \left( g - P^s_{B_j} g \right) \chi_{B_j} \right\|_{L^{\vec{p}'}(\mathbb{R}^n)} \|a_j\|_{L^{\vec{p}}(\mathbb{R}^n)} \leqslant \|g\|_{\mathcal{CL}_{\vec{p}',\alpha,s}(\mathbb{R}^n)},$$

which implies $\mathcal{CL}_{\vec{p}',\alpha,s}(\mathbb{R}^n) \subseteq \left( H\dot{K}^{\alpha,q}_{\vec{p}}(\mathbb{R}^n) \right)^*$.

(ii) $\left( H\dot{K}^{\alpha,q}_{\vec{p}}(\mathbb{R}^n) \right)^* \subseteq \mathcal{CL}_{\vec{p}',\alpha,s}(\mathbb{R}^n)$.

For this part, it suffices to check that for any $L \in \left( H\dot{K}^{\alpha,q}_{\vec{p}}(\mathbb{R}^n) \right)^*$, there exist some $g \in \mathcal{CL}_{\vec{p}',\alpha,s}(\mathbb{R}^n)$ such that for all $f = \sum_{i=-\infty}^{\infty} \lambda_i a_i \in \left( H\dot{K}^{\alpha,q}_{\vec{p}}(\mathbb{R}^n) \right)^*$,

$$L(f) = \sum_{i=-\infty}^{\infty} \lambda_i \int_{\mathbb{R}^n} a_i(x) g(x) dx.$$

For any measurable set $E$, we define a function space

$$L^{\vec{p}}_0(E) = \left\{ f \in L^{\vec{p}}(E) : \int_E f(x) \cdot x^\beta dx = 0 \right\} \quad \text{when} \quad |\beta| \leqslant s.$$

Obviously, the function $|B_k|^{-\frac{\alpha}{n}} \|f\|^{-1}_{L^{\vec{p}}(B_k)} \cdot f$ is a central $(\alpha, \vec{p})$ atom supported in $B_k$.

Since $L$ is a bounded linear functional on $H\dot{K}^{\alpha,q}_{\vec{p}}(\mathbb{R}^n)$, for any $f \in L^{\vec{p}}_0(B_k)$, we have

$$|L(f)| \leqslant \|L\| \cdot \|f\|_{H\dot{K}^{\alpha,q}_{\vec{p}}(\mathbb{R}^n)} \leqslant |B_k|^{\frac{\alpha}{n}} \|f\|_{L^{\vec{p}}(B_k)}.$$

By Hahn-Banach theorem, there exists a function $g_k \in L^{\vec{p}'}(\mathbb{R}^n)$ such that

$$L(f) = \int_{B_k} f(x) g_k(x) dx.$$

Therefore, if we can find a function $g$ such that for any $l \in \mathbb{Z}$ and $f \in L^{\vec{p}}_0(B_l)$,

$$L(f) = \int_{B_l} f(x) g(x) dx,$$

then $g(x) \in \left( H\dot{K}^{\alpha,q}_{\vec{p}}(\mathbb{R}^n) \right)^*$, and in addition,

$$L(f) = \sum_{i=-\infty}^{\infty} \lambda_i \int_{\mathbb{R}^n} a_i(x) g(x) dx.$$

Next we will consider that existence of $g$. For $f \in L^{\vec{p}}_0(B_1)$, there exists $g_1$ such that

$$L(f) = \int_{B_1} f(x) g_1(x) dx.$$

For $f \in L_0^{\vec{p}}(B_2)$, there exist $g_2$ satisfying

$$L(f) = \int_{B_2} f(x) g_2(x) dx.$$

Noting that $B_1 \subseteq B_2$, for $\forall f \in L_0^{\vec{p}}(B_1)$, we obtain

$$\int_{B_1} f(x) \left(g_1(x) - g_2(x)\right) dx = 0.$$

If we replace $f$ by any $h(x) - P_{B_1}^s h(x)$ where $h \in L^{\vec{p}}(\mathbb{R}^n)$, then

$$\int_{B_1} \left(h(x) - P_{B_1}^s h(x)\right) \left(g_1(x) - g_2(x)\right) dx = 0,$$

Replacing $g_1(x) - g_2(x)$ by $(g_1(x) - g_2(x)) - P_{B_1}^s (g_1 - g_2)(x)$, we have

$$\int_{B_1} \left(h(x) - P_{B_1}^s h(x)\right) \left(g_1(x) - g_2(x) - P_{B_1}^s (g_1 - g_2)(x)\right) dx = 0.$$

As a consequence, for all $x \in B_1$, $g_1 = g_2 + P_{B_1}^s (g_1 - g_2)$.

Denote

$$g(x) = \begin{cases} g_1(x) & x \in B_1, \\ g_2 + P_{B_1}^s (g_1 - g_2)(x) & x \in B_2 \end{cases}$$

For $g$ defined above, the following equality hold:

$$L(f) = \int_{B_i} f(x) g(x) dx \quad \text{for all } f \in L_0^{\vec{p}}(B_i).$$

By repeating the above process, we can get a function $h$ such that

$$L(f) = \int_{B_i} f(x) g(x) dx \quad (i = 1, 2 \cdots).$$

Next, we only need to prove $g \in \mathcal{CL}_{\vec{p}, \alpha, s}(\mathbb{R}^n)$. In fact,

$$\left\| (g - P_{B_k}^s g) \chi_{B_k} \right\|_{L^{\vec{p}'}} = \sup_{\|f\|_{L^{\vec{p}}}=1} \left| \int_{B_k} \left(g(x) - P_{B_k}^s g(x)\right) f(x) dx \right|$$

$$= \sup_{\|f\|_{L^{\vec{p}}}=1} \left| \int_{B_k} g(x) f(x) dx \right|$$

$$\lesssim \sup_{\|f\|_{L^{\vec{p}}}=1} |L(f)|$$

$$\lesssim \|L\| \|f\|_{L^{\vec{p}}} |B|^{\frac{\alpha}{n}}.$$

Therefore,

$$\left\| (g - P_{B_k}^s g) \right\|_{L^{\vec{p}}} |B|^{-\frac{\alpha}{n}} \lesssim \|L\|.$$

This proof is complete. □

**Competing interests**

The authors declare that they have no competing interests.

**Funding**

The research was supported by Natural Science Foundation of China (Grant No. 12061069) and the Natural Science Foundation of Henan Province(Grant No. 202300410338 ).

**Authors contributions**

All authors contributed equality and significantly in writing this paper. All authors read and approved the final manuscript.

**Acknowledgments**

The authors would like to express their thanks to the referees for valuable advice regarding previous version of this paper.

**Authors detials**

Yichun Zhao and Jiang Zhou*, 13565901944@163.com and zhoujiang@xju.edu.cn, College of Mathematics and System Science, Xinjiang University, Urumqi, 830046, P.R China.

Mingquan Wei, weimingquan11@mails.ucas.ac.cn, School of Mathematics and Statistics, Xinyang Normal University, Xinyang, 464000, P.R China.

Yichun Zhao and Jiang Zhou

College of Mathematics and System Sciences

Xinjiang University

Urumqi 830046

E-mail : `13565901944@163.com` (Yichun Zhao); `zhoujiang@xju.edu.cn` (Jiang Zhou)

Mingquan Wei

School of Mathematics and Stastics

Xinyang Normal University

Xinyang, 464000

E-mail : `weimingquan11@mails.ucas.ac.cn` (Mingquan Wei)